\pdfoutput=1
\RequirePackage{ifpdf}
\ifpdf 
\documentclass[pdftex]{sigma}
\else
\documentclass{sigma}
\fi

\usepackage{mathrsfs}
\usepackage{enumitem}
\usepackage{mathtools}

\DeclareMathOperator{\tr}{tr}
\DeclareMathOperator{\dvol}{dvol}
\DeclareMathOperator{\Ric}{Ric}
\DeclareMathOperator{\Rm}{Rm}

\newcommand{\cg}{\widetilde{g}}
\newcommand{\cR}{\widetilde{R}}
\newcommand{\cf}{\widetilde{f}}
\newcommand{\cmG}{\widetilde{\mathcal{G}}}
\newcommand{\cRic}{\widetilde{\Ric}}

\newcommand{\lp}{\langle}
\newcommand{\rp}{\rangle}
\newcommand{\lv}{\lvert}
\newcommand{\rv}{\rvert}

\newcommand{\mF}{\mathcal{F}}
\newcommand{\mG}{\mathcal{G}}
\newcommand{\mW}{\mathcal{W}}

\newcommand{\bN}{\mathbb{N}}
\newcommand{\bR}{\mathbb{R}}
\newcommand{\del}{\delta_{\phi}}
\newcommand{\ric}{\Ric_{\phi}^m}
\newcommand{\psmms}{(M_+,g_+,f_+,m,\mu)}

\numberwithin{equation}{section}

\newtheorem{Theorem}{Theorem}[section]

\newtheorem{Lemma}[Theorem]{Lemma}
\newtheorem{Proposition}[Theorem]{Proposition}
 { \theoremstyle{definition}
\newtheorem{Definition}[Theorem]{Definition}

\newtheorem{Remark}[Theorem]{Remark} }

\begin{document}
\allowdisplaybreaks

\newcommand{\arXivNumber}{2202.11153}

\renewcommand{\PaperNumber}{086}

\FirstPageHeading

\ShortArticleName{The Weighted Ambient Metric}

\ArticleName{The Weighted Ambient Metric}

\Author{Jeffrey CASE and Ayush KHAITAN}

\AuthorNameForHeading{J.S.~Case and A.~Khaitan}

\Address{Department of Mathematics, Penn State University, PA, USA}
\Email{\href{mailto:jscase@psu.edu}{jscase@psu.edu}, \href{mailto:auk480@psu.edu}{auk480@psu.edu}}
\URLaddress{\url{http://www.personal.psu.edu/jqc5026/}}

\ArticleDates{Received March 01, 2022, in final form October 25, 2022; Published online November 02, 2022}

\Abstract{We prove the existence and uniqueness of weighted ambient metrics and weighted Poincar\'e metrics for smooth metric measure spaces.}


\Classification{53A31; 53A55; 31C12}

\section{Introduction}

The Fefferman--Graham ambient space~\cite{FeffermanGraham2012} is a formally Ricci flat space canonically associated to a given conformal manifold.
Among its many applications are the classification of local scalar conformal invariants~\cite{BaileyEastwoodGraham1994,Gover2001} and the construction of a family of conformally covariant operators, called GJMS operators, with leading-order term a power of the Laplacian~\cite{GJMS1992}.
One can also use the ambient space to construct an asymptotically hyperbolic, formally Einstein space with conformal boundary a given conformal manifold.
The resulting space, called a Poincar\'e space, likewise has many applications, among them a proof that the GJMS operators are formally self-adjoint~\cite{GrahamZworski2003} and the construction and study of variational scalar invariants generalizing the $\sigma_2$-curvatures~\cite{ChangFang2008,ChangFangGraham2012,Graham2009}.

A \emph{smooth metric measure space} is a five-tuple $(M^d,g,f,m,\mu)$ formed from a Riemannian manifold $(M^d,g)$, a positive function $f\in C^\infty(M)$, a dimensional parameter $m\in[0,\infty]$, and a curvature parameter $\mu\in\bR$. Smooth metric measure spaces arise in many ways, including as (possibly collapsed) limits of sequences of Riemannian manifolds~\cite{CheegerColding2000a}, as smooth manifolds satisfying curvature-dimension inequalities~\cite{BakryEmery1985,Wei_Wylie}, as the geometric framework~\cite{CaseChang2013} for studying curved analogues of the Caffarelli--Silvestre extension for defining the fractional Laplacian~\cite{CaffarelliSilvestre2007}, and, in the limiting case $m=\infty$, as a geometric framework for the realization of the Ricci flow as a gradient flow~\cite{Perelman1}.

A \emph{weighted invariant} is tensor-valued function on smooth metric measure spaces such that $T(N,\Phi^\ast g,\Phi^\ast f,m,\mu) \!=\! \Phi^\ast(T(M,g,f,m,\mu))$ for any smooth metric measure space $(M^d,g,f,m,\mu)$ and any diffeomorphism $\Phi \colon N \mapsto M$.
For example, the \emph{Bakry--\'Emery Ricci curvature} of $(M^d,g,f,m,\mu)$ is
\begin{equation*}
 \Ric_\phi^m := \Ric - \frac{m}{f}\nabla^2f
\end{equation*}
and the \emph{weighted scalar curvature} of $(M^d,g,f,m,\mu)$ is
\begin{equation*}
 R_\phi^m := R - \frac{2m}{f}\Delta f - \frac{m(m-1)}{f^2}\bigl( \lv\nabla f\rv^2 - \mu \bigr) .
\end{equation*}
A smooth metric measure space is often studied as an abstract analogue of the base of the warped product $M^d \times_f F^m(\mu)$ of $(M^d,g)$ with an $m$-dimensional spaceform $F^m(\mu)$ of constant sectional curvature $\mu$, in the sense that when $m$ is a nonnegative integer, most weighted invariants are equivalent to Riemannian invariants of the warped product.
For example, if $m \in \bN_0$, then $\Ric_\phi^m$ is the restriction of the Ricci tensor of $M \times_f F^m(\mu)$ to horizontal (i.e., tangent to $M$) vectors and $R_\phi^m$ is the scalar curvature of $M \times_f F^m(\mu)$.

We say that two smooth metric measure spaces $(M^d,g_i,f_i,m,\mu)$, $i\in\{1,2\}$, are \emph{pointwise conformally equivalent} if there is a function $u\in C^\infty(M)$ such that $g_2={\rm e}^{2u}g_1$ and $f_2={\rm e}^uf_1$.
This notion has many applications. For example, the curved version of the Caffarelli--Silvestre extension~\cite{CaffarelliSilvestre2007,Case2019tr, CaseChang2013} identifies a conformally covariant operator with leading-order term a fractional power of the Laplacian~\cite{GrahamZworski2003} as a generalized Dirichlet-to-Neumann operator associated to a weighted analogue of the GJMS operators~\cite{Case2011t, CaseChang2013}, and there is a conformally invariant analogue of the Yamabe functional on smooth metric measure spaces~\cite{Case2013y} which interpolates between the usual Yamabe functional~\cite{Yamabe1960} and Perelman's $\mF$- and $\mW$-functionals~\cite{Perelman1}.
This latter perspective leads to an ad hoc construction~\cite{Case2014sd,Case2016v, Case2014s} of some fully nonlinear analogues of Perelman's functionals.

Given the variety of applications of the Fefferman--Graham ambient space and of the conformal geometry of smooth metric measure spaces, it is natural to ask whether there is a canonical weighted ambient space associated to a smooth metric measure space. We show that this is the case.
Roughly speaking, we show that if $(M^d,[g,f],m,\mu)$, $d\geq 3$ and $m<\infty$, is a conformal class of smooth metric measure spaces, then there is a unique smooth metric measure space $(\mathbb{R}_+\times M\times (-\epsilon,\epsilon),\widetilde{g},\widetilde{f},m,\mu)$ of Lorentzian signature such that:
\begin{enumerate}[label=$(\roman*)$]\itemsep=0pt
\item if $d+m$ is not an even integer, then $\widetilde{\mathrm{Ric}}_\phi^m,\widetilde{F_\phi^m}=O(\rho^\infty)$, where $\rho$ is the coordinate on $(-\epsilon,\epsilon)$ and
 \[
 \widetilde{F_\phi^m}:=\widetilde{f}\widetilde{\Delta}\widetilde{f}+(m-1)\bigl(\big|\widetilde{\nabla} \widetilde{f}\big|^2-\mu\bigr);
 \]

\item if $d+m$ is an even integer, then
\begin{gather*}
 \widetilde{\mathrm{Ric}}_\phi^m,\widetilde{F_\phi^m}=O\bigl(\rho^{\frac{d+m}{2}-1}\bigr),
 \qquad
 g^{ij}\bigl(\widetilde{\mathrm{Ric}}_\phi^m\bigr)_{ij}-{m}{f}^{-2} \widetilde{F_\phi^m}=O\bigl(\rho^{\frac{d+m}{2}}\bigr),
\end{gather*}
where $i$, $j$ denote the coordinates on $M$.
\end{enumerate}
Note that $\widetilde{R_\phi^m}=\widetilde{g}^{IJ}\big(\widetilde{\mathrm{Ric}}_\phi^m\big)_{IJ}-m f^{-2}\widetilde{F_\phi^m}$, where $I$, $J$ denote the coordinates on $\widetilde{\mathcal{G}}:=\mathbb{R}_+\times M\times (-\epsilon,\epsilon)$. The improved orders of vanishing in the case that $d+m$ is even can be interpreted as the statement that, modulo the $\rho$ components of $\widetilde{g}^{IJ}\big(\widetilde{\mathrm{Ric}}_\phi^m\big)_{IJ}$,
the weighted scalar curvature~\smash{$\cR_\phi^m$} vanishes to one higher order in $\rho$ than the \smash{Bakry--\'Emery} Ricci tensor \smash{$\cRic_\phi^m$}. Our arguments carry through verbatim if $g$ has signature $(p,q)$, in which case $\widetilde{g}$ has signature $(p+1,q+1)$.

A precise statement of the above result, using the terminology of Section~\ref{sec:fg}, is as follows:

\begin{Theorem}
 \label{ambient-statement}
Let $(M^d, [g, f], m, \mu)$, $d\geq 3$ and $m<\infty$, be a conformal class of smooth metric measure spaces. Then
\begin{enumerate}[label=$(\roman*)$]\itemsep=0pt
\item there exists a unique, up to ambient-equivalence, weighted ambient space $(\widetilde{\mathcal{G}}, \widetilde{g},\widetilde{f},m,\mu)$ for $(M^d,[g,f],m,\mu)$; and

\item if $d+m$ is an even integer, the weighted obstruction tensor,
\begin{align}
&\mathcal{O}_{ij}:=c_{d+m}\partial_\rho^{\frac{d+m}{2}-1}\big|_{\rho=0} \big(\widetilde{{\rm Ric}_\phi^m}\big)_{ij},\nonumber
\\
&c_{d+m}=(-2)^{(d+m)/2-1}\frac{((d+m)/2-1)!}{d+m-2},\label{obstruction-definition}
\end{align}
is a local conformal invariant of $(M^d,[g,f],m\mu)$ which vanishes if and only if there is a~weighted ambient space $(\widetilde{g},\widetilde{f})$ such that $\widetilde{\mathrm{Ric}_\phi^m},\widetilde{F_\phi^m}=O(\rho^{\infty})$.
\end{enumerate}
\end{Theorem}

When $m=0$, Theorem~\ref{ambient-statement} recovers the Fefferman--Graham ambient space~\cite{FeffermanGraham2012}.
When~$m$ is a positive integer, the existence and uniqueness statement of Theorem~\ref{ambient-statement} implies that $\cmG \times_{\cf} F^m(\mu)$ is the Fefferman--Graham ambient space of $M \times_f F^m(\mu)$. See Theorem~\ref{obstruction-theorem} for additional properties of the weighted obstruction tensor.

We do not presently know if there is an analogue of Theorem~\ref{ambient-statement} in the case $m=\infty$.

Further expansions involving fractional powers of $\rho$ when $d+m\notin 2\mathbb{N}$ or $\log$ terms when $d+m\in\mathbb{N}$ are also possible~(cf.~\cite{FeffermanGraham2012}).
We have not developed this because it is unnecessary for our intended applications to the construction of local invariants.

From Section~\ref{section-4.1} onwards, the metric measure structure $(\widetilde{g},\widetilde{f})$ will be assumed to be of the form
\begin{align}
\label{weighted-straight-normal-form}
&\widetilde{g}=2\rho \,\mathrm{d}t^2+2t \,{\rm d}\rho \,\mathrm{d}t+ t^2 g_\rho,\qquad
\widetilde{f}= tf_\rho,
\end{align}
where $(g_\rho,f_\rho)$ is a one-parameter family of metric measure structures on $M$.
Metrics of the form equation~\eqref{weighted-straight-normal-form} are called \emph{straight} and \emph{normal}. See Sections~\ref{sec:fg} and~\ref{formal-theory} for more details.

It is in general laborious to compute the asymptotic expansions of $(\widetilde{g},\widetilde{f})$.
However, one readily computes that, in the case of a straight and normal weighted ambient space,
\begin{gather*}
 g_\rho = g + 2\rho P_\phi^m+O\big(\rho^2\big), \qquad
 f_\rho = f + \frac{f}{m}\rho Y_\phi^m +O\big(\rho^2\big),
\end{gather*}
where
\begin{equation*}
 P_\phi^m := \frac{1}{d+m-2}\bigl( \Ric_\phi^m - J_\phi^m g\bigr)
\end{equation*}
is the \emph{weighted Schouten tensor},
\begin{equation*}
 J_\phi^m := \frac{1}{2(d+m-1)}R_\phi^m
\end{equation*}
is the \emph{weighted Schouten scalar}, and $Y_\phi^m := J_\phi^m - \tr_g P_\phi^m$.

We can also explicitly identify the weighted ambient space in three special cases.

\begin{Theorem}
 \label{examples-theorem}
Let $(M^d,g,f,m,\mu)$, $d\geq 3$ and $m<\infty$, be a smooth metric measure space. Set $\widetilde{\mathcal{G}}=\mathbb{R}_+\times M^d\times (-\epsilon,\epsilon)$ and consider the smooth metric measure space $(\widetilde{\mathcal{G}},\widetilde{g},\widetilde{f},m,\mu)$, where $(\widetilde{g},\widetilde{f})$ is of the form of equation~\eqref{weighted-straight-normal-form} and
\begin{gather*}
(g_{\rho})_{i j}=g_{i j}+2\rho\big(P_{\phi}^{m}\big)_{i j} +\rho^{2}\big(P_{\phi}^{m}\big)_{i k}\big(P_{\phi}^{m}\big)_{j}^{k},
\\
f_{\rho}=f+ \frac{f}{m}\rho Y_{\phi}^{m}.
\end{gather*}
Then $(\widetilde{\mathcal{G}},\widetilde{g},\widetilde{f},m,\mu)$ is a weighted ambient space if
\begin{enumerate}[label=$(\roman*)$]\itemsep=0pt
\item \label{Theorem 1.2-1} there is a constant $\lambda\in\mathbb{R}$ such that
\[
F_{\phi}^{m}=- 2(d+m-1)\lambda f^{2}\hspace{5mm} \text { and }\hspace{5mm} \mathrm{Ric}_{\phi}^{m}= 2(d+m-1)\lambda g;
\]

\item \label{Theorem 1.2-2} $Rm=P_\phi^m\wedge g$ and ${\rm d}P_\phi^m=0$; or

\item \label{Theorem 1.2-3} $f=1$ and $\mathrm{Ric}=-(d-1)\mu g$.
\end{enumerate}
\end{Theorem}

The first case of Theorem~\ref{examples-theorem} is the case of a quasi-Einstein manifold~\cite{Case_Shu_Wei};
the second case is the case of a smooth metric measure space which is locally conformally flat in the weighted sense~\cite{Case2014s};
and the third case is a weighted interpretation of the explicit ambient metric of Gover and Leitner~\cite{GoverLeitner2009}.

As in Riemannian geometry, the weighted ambient space gives rise to a canonical weighted Poincar\'e space associated to a given smooth metric measure space.

\begin{Theorem}
\label{poincare-statement}
Let $(M^{d}, [g, f], m, \mu)$, $d\geq 3$ and $m<\infty$, be a conformal class of smooth metric measure spaces. Then there exists an even weighted Poincar\'e space for it. Moreover, if $\bigl(g_{+}^{1}, f_{+}^{1}\bigr)$ and $\bigl(g_{+}^{2}, f_{+}^{2}\bigr)$ are two even weighted Poincar\'e spaces for $\bigl(M^{d}, g, f, m, \mu\bigr)$ defined on $\bigl(M_{+}^{1}\bigr)^{\circ}$ and~$\bigl(M_{+}^{2}\bigr)^{\circ}$, respectively, then there are open sets $U_{1} \subset M_{+}^{1}$ and $U_{2} \subset M_{+}^{2}$ containing $M \times\{0\}$ and an even diffeomorphism $\phi\colon U_{1} \rightarrow$ $U_{2}$ such that $\phi|_{M \times\{0\}}$ is the identity map and
\begin{enumerate}[label=$(\roman*)$]\itemsep=0pt
\item if $m+d$ is not an even integer, then
\[
\bigl(g_{+}^{1}-\phi^{*} g_{+}^{2}, -2\big(f_{+}^{1}-\phi^{*} f_{+}^{2}\big)\bigr)=O(r^{\infty});
\]

\item if $d+m$ is an even integer, then
\[
\bigl(g_{+}^{1}-\phi^{*} g_{+}^{2}, -2\big(f_{+}^{1}-\phi^{*} f_{+}^{2}\big)\bigr)=O_{\alpha \beta}^{1,+}\big(r^{d+m-2}\big).
\]
\end{enumerate}
\end{Theorem}
See Sections~\ref{sec:fg} and~\ref{sec:poincare} for the relevant definitions.

In a separate article \cite{Khaitan2022}, Khaitan used Theorems~\ref{ambient-statement}, \ref{examples-theorem} and~\ref{poincare-statement} to give a rigorous construction of the weighted GJMS operators of all orders up to the obstruction, and to prove that they are formally self-adjoint. He also gave explicit formulas for the weighted GJMS operators for smooth metric measure spaces satisfying the first or last conditions of Theorem~\ref{examples-theorem}.
This makes rigorous the factorization of the weighted GJMS operators established by Case and Chang~\cite{CaseChang2013} in the latter case by formally arguing via warped products.

This article is organized as follows: In Section~\ref{sec:bg.tex}, we recall some weighted invariants of smooth metric measure spaces. In Section~\ref{sec:fg}, we define weighted ambient spaces and discuss weighted ambient-equivalence. In Section~\ref{formal-theory}, we prove Theorem~\ref{ambient-statement}. In Section~\ref{sec:poincare}, we define weighted Poincar\'e spaces and prove Theorem~\ref{poincare-statement}. In Section~\ref{sec:explicit-examples}, we prove Theorem~\ref{examples-theorem}.

\section{Smooth metric measure spaces}
\label{sec:bg.tex}

We recall some weighted invariants relevant to the geometry of smooth metric measure spaces~\cite{Case2014s}. Let $(M^d,g,f,m,\mu)$ be a smooth metric measure space. We discuss here only the case $m<\infty$ due to the presence of this restriction to our main results. This data determines a volume element
\[
{\rm d}v_\phi := f^m\,\dvol_g .
\]
When $m>0$, we set $\phi:=-m\ln f$, so that
\[
{\rm d}v_\phi = {\rm e}^{-\phi}\,\dvol_g .
\]

In terms of $\phi$, it holds that
\begin{align*}
&\Ric_\phi^m= \Ric + \nabla^2\phi - \frac{1}{m} {\rm d}\phi \otimes {\rm d}\phi,
\\
&R_\phi^m= R + 2\Delta\phi - \frac{m+1}{m}\lv\nabla\phi\rv^2 + m(m-1)\mu {\rm e}^{2\phi/m} .
\end{align*}
Recall that \emph{weighted Schouten tensor $P_\phi^m$} and the \emph{weighted Schouten scalar $J_\phi^m$} of $(M^d,g,f,m,\mu)$ are
\begin{gather*}
 P_\phi^m := \frac{1}{d+m-2}\big(\Ric_\phi^m - J_\phi^m g\big), \qquad
 J_\phi^m := \frac{1}{2(d+m-1)}R_\phi^m .
\end{gather*}
We set
\[
F_\phi^m:= f\Delta f+(m-1)\big(|\nabla f|^2-\mu\big).
\]
Observe that
\begin{equation}\label{fphim-equation}
F_\phi^m=\frac{f^2}{m}\bigl[(d+m-2)\mathrm{tr}_g P_\phi^m -(d+2m-2)J_\phi^m\bigr].
\end{equation}

A \emph{weighted Einstein manifold} is a smooth metric measure space $(M^d,g,f,m,\mu)$ such that
\[
P_\phi^m = \lambda g
\]
for some constant $\lambda\in\bR$.
A weighted Einstein manifold is \emph{quasi-Einstein}~\cite{Case_Shu_Wei} if additionally
\[
J_\phi^m=(d+m)\lambda.
\]
It follows from equation~\eqref{fphim-equation} that this is equivalent to
\[
\Ric_\phi^m=2(d+m-1)\lambda g,\qquad F_\phi^m=-2(d+m-1)\lambda f^2.
\]

We define the \emph{weighted divergence} $\delta_\phi T$ of a tensor field $T \in \Gamma(\otimes^k T^*M)$ by
\[
(\delta_\phi T)(X_1,\dots,X_k):=\sum\limits_{i=1}^d \nabla_{e_i} T(e_i,X_1,\dots,X_k)-T(\nabla\phi,X_1,\dots,X_k),
\]
where $X_1,\dots,X_k\in T_p M$ and $\{e_i\}$ is an orthonormal basis for $T_p M$.

The \emph{weighted Weyl tensor $A_\phi^m$} and the \emph{weighted Cotton tensor ${\rm d}P_\phi^m$} of a smooth metric measure space $(M^d,g,f,m,\mu)$ are
\begin{gather*}
 A_\phi^m(x,y,z,w) := \big(\Rm - P_\phi^m \wedge g\big)(x,y,z,w) ,
 \\
 {\rm d}P_\phi^m(x,y,z) := \nabla_x P_\phi^m(y,z) - \nabla_y P_\phi^m(x,z) ,
\end{gather*}
where $h\wedge k$ denotes the Kulkarni--Nomizu product
\[
(h\wedge k)(x,y,z,w) := h(x,z)k(y,w) + h(y,w)k(x,z) - h(x,w)k(y,z) - h(y,z)k(x,w) .
\]
The \emph{weighted Bach tensor} of $(M^d,g,f,m,\mu)$ is
\[
B_\phi^m := \delta_\phi {\rm d}P_\phi^m - \frac{1}{m}\tr {\rm d}P_\phi^m \otimes {\rm d}\phi + \biggl\lp A_\phi^m, P_\phi^m - \frac{Y_\phi^m}{m}g\biggr\rp ,
\]
where $Y_\phi^m := J_\phi^m - \tr P_\phi^m$ and the contractions are
\begin{gather*}
 (\tr {\rm d}P_\phi^m)(x) := \sum_{i=1}^d {\rm d}P_\phi^m(e_i,x,e_i),
 \\
 \biggl\lp A_\phi^m, P_\phi^m - \frac{Y_\phi^m}{m}g\biggr\rp(x,y) := \sum_{i,j=1}^d A_\phi^m(e_i,x,e_j,y)\biggl(P_\phi^m - \frac{Y_\phi^m}{m}g\biggr)(e_i,e_j).
\end{gather*}

The \emph{weighted Bianchi identity}~\cite{Case2014sd} is
\begin{equation}
\label{intro-weighted-bianchi}
 \del \ric-\frac{1}{2}{\rm d}P_\phi^m-\frac{1}{f^2} F_\phi^m {\rm d}\phi=0.
\end{equation}

{\samepage
A smooth metric measure space $(M^d,g,f,m,\mu)$ is \emph{locally conformally flat in the weighted sense}~\cite{Case2014s} if
\begin{enumerate}[label=$(\roman*)$]\itemsep=0pt
 \item $d+m\in\{1,2\}$;
 \item $d+m=3$ and ${\rm d}P_\phi^m=0$; or
 \item $d+m\not\in\{1,2,3\}$ and $A_\phi^m=0$.
\end{enumerate}

}

Note that~\cite{Case2014s} if $d+m\geq 3$, then $(M^d,g,f,m,\mu)$ is locally conformally flat in the weighted sense if and only if $A_\phi^m,{\rm d}P_\phi^m=0$.

We require some formulas satisfied by smooth metric measure spaces which are locally conformally flat in the weighted sense.
\begin{Lemma}[{\cite[Lemma~3.2]{Case2014s}}]
\label{locally-conformally-flat-implications}
Let $(M^{d}, g, v, m, \mu)$ be a smooth metric measure space that is locally conformally flat in the weighted sense. Then
\begin{subequations}
\begin{gather}
0=P_{\phi}^{m}(\nabla \phi)+{\rm d}Y_{\phi}^{m}-\frac{1}{m} Y_{\phi}^{m} {\rm d}\phi,
\\
\label{second-condition-locally-conformally-flat}
0=m P_{\phi}^{m}-\nabla^{2} \phi+\frac{1}{m} {\rm d}\phi \otimes {\rm d}\phi+ Y_{\phi}^{m} g,
\\
0={\rm d} P_{\phi}^{m}.
\end{gather}
\end{subequations}
\end{Lemma}

\section{Weighted ambient spaces}
\label{sec:fg}

Let $(M^d,g,f,m,\mu)$ be a smooth metric measure space with $d\geq 3$ and $m<\infty$. Denote by $\mathcal{E}$ the trivial line bundle over $M$. Let $\mathcal{G}$ be the ray subbundle of $S^2 T^*M\oplus \mathcal{E}$ consisting of all triples $(h,u,x)$ such that $h=s^2 g_x$ and $u=sf(x)$ for some $s\in \bR_+$ and $x\in M$. We define the dilation $\delta_s\colon \mG \mapsto \mG$ by $(h,u,x)\mapsto (s^2 h,su,x)$. Let $T:=\frac{{\rm d}}{{\rm d}s}\delta_s\big|_{s=1}$ denote the infinitesimal generator of dilations. Also, let $\pi \colon (h,u,x) \mapsto x$ be the projection from $\mG$ to $M$. There is a canonical metric measure structure $(\boldsymbol{g},\boldsymbol{f})$ on $\mG$ defined by $\boldsymbol{g}(X,Y)=h(\pi_*X,\pi_*Y)$ and $\boldsymbol{f}(h,u,x)=u(x)$.

Given $(g,f)$, define a coordinate chart $\mG\mapsto \bR_+\times M$, by $(t^2g_x,tf(x),x)\mapsto (t,x)$. Then $T=t\partial_t$.

Consider the embedding $\iota \colon \mG\ \hookrightarrow \mG\times\bR,(t,x)\mapsto (t,x,0)$. Each point in $\mG\times\bR$ can be written as $(t,x,\rho)$. We extend the dilation to $\mathcal{G}\times \bR$ as $\delta_s (t,x,\rho)= (st,x,\rho)$.

We now define weighted pre-ambient spaces and special cases thereof.

\begin{Definition}
$(\cmG,\cg,\cf,m,\mu)$ is called a \emph{weighted pre-ambient space} if
\begin{enumerate}[label=$(\roman*)$]\itemsep=0pt
 \item $\cmG$ is a dilation-invariant neighborhood of $\mG\times \{0\}$;
 \item $\cg$ is a smooth metric of signature $(d+1,1)$ and $\widetilde{f}$ is a smooth function on $\cmG$;
 \item $\delta_s^*\widetilde{g}=s^2 \widetilde{g}$ and $\delta_s^*\widetilde{f}=s\widetilde{f}$; and
 \item $(\iota^*\widetilde{g},\iota^*\widetilde{f})=(\boldsymbol{g},\boldsymbol{f})$.
\end{enumerate}
\end{Definition}

\begin{Remark}
Our results can be extended to pseudo-Riemannian metrics in the usual way (cf.~\cite{FeffermanGraham2012}).
\end{Remark}

\begin{Definition}
A weighted pre-ambient space $(\cmG,\cg,\cf,m,\mu)$ is said to be in \emph{normal form} relative to a metric measure structure $(g,f)$ if
\begin{enumerate}[label=$(\roman*)$]\itemsep=0pt
 \item for each fixed $z\in\mG$, the set of $\rho\in\bR$ such that $(z,\rho)\in \cmG$ is an open interval $I_z$ containing~$0$;
 \item for each $z\in\mG$, the curve on $I_z$ defined as $\rho\mapsto(z,\rho)$ is a geodesic in $\cmG$; and
 \item on $\mG\times \{0\}$, it holds that $\cg=\boldsymbol{g} + 2t\,{\rm d}\rho\,\mathrm{d}t$.
\end{enumerate}
\end{Definition}

\begin{Definition}
A weighted pre-ambient space $(\cmG,\cg,\cf,m,\mu)$ for $(M^d,[g,f],m,\mu)$ is be said to be \textit{straight} if any of the following equivalent properties hold:
\begin{enumerate}[label=$(\roman*)$]\itemsep=0pt
 \item for each $p\in \widetilde{\mathcal{G}}$, the dilation orbit $s\mapsto \delta_s p$ is a geodesic for $\widetilde{g}$;

 \item $\widetilde{g}(2T,\cdot)={\rm d}(\widetilde{g}(T,T))$; or

\item the infinitesimal dilation field $T$ satisfies $\widetilde{\nabla}T=\mathrm{Id}$.
\end{enumerate}
\end{Definition}

Given a point $(t,x,\rho)\in \mathcal{G}\times \bR$, we refer to the $t$ coordinate as the $0$ coordinate, and the $\rho$ coordinate as the $\infty$ coordinate. Coordinates in $M$ are denoted with lowercase latin characters $(i,j,k,\dots)$.

The metric of a normal weighted pre-ambient space takes a special form.

\begin{Lemma}
\label{straight-metric-lemma}
Let $(\cmG,\cg,\cf,m,\mu)$ be a weighted pre-ambient space such that for each $z\in\mG$, the set of all $\rho\in\bR$ such that $(z,\rho)\in\cmG$ is an open interval containing $0$. Then $(\cg,\cf)$ is a \textit{normal} metric measure structure if and only if $\cg_{0\infty}=t$, $\cg_{\infty i}=0$ and $\cg_{\infty\infty}=0$.
\end{Lemma}

\begin{proof}
A normal metric measure structure has
\[
\widetilde{g}_{\infty\infty}\big|_{\rho=0}=\widetilde{g}_{i\infty}\big|_{\rho=0}=0 \qquad\text{and}\qquad \widetilde{g}_{0\infty}\big|_{\rho=0}=t.
\]
Moreover, it also has geodesic $\rho$-lines. The $\rho$-lines are geodesics if and only if the Christoffel symbols $\widetilde{\Gamma}_{\infty \infty I}$, $I\in \{0,i,\infty\}$, vanish. Taking $I=\infty$ gives $\partial_{\rho} \widetilde{g}_{\infty \infty}=0$, and hence $\widetilde{g}_{\infty \infty}=0 .$ Now taking $I=i$ and $I=0$ gives $\widetilde{g}_{\infty i}=0$ and $\widetilde{g}_{\infty 0}=t$.

The converse follows similarly.
\end{proof}
Defining a weighted ambient space required some additional notation.

\begin{Definition}
Let $(\widetilde{S}_{IJ},\widetilde{h})$ be a pair of a symmetric 2-tensor $\widetilde{S}_{IJ}$ and a smooth function $\widetilde{h}$ on an open neighborhood of $\mathcal{G}\times \{0\}\subset\mathcal{G}\times \bR$. For $k,m\geq 0$, we write $(\widetilde{S}_{IJ},\widetilde{h})=O_{IJ}^{k,+}(\rho^m)$ if
\begin{enumerate}[label=$(\roman*)$]\itemsep=0pt
\item $\big(\widetilde{S}_{IJ},\widetilde{h}\big)=O(\rho^{m})$;
\item $\widetilde{S}_{00},\widetilde{S}_{0i}=O(\rho^{m+1})$; and
\item $m f^{-k}\widetilde{h}-g^{ij}\widetilde{S}_{ij}= O(\rho^{m+1})$.
\end{enumerate}
\end{Definition}

We now define a weighted ambient space.

\begin{Definition}
A \emph{weighted ambient space} for $(M^d,g,f,m,\mu)$ is a weighted pre-ambient space $(\cmG,\cg,\cf,m,\mu)$ such that
\begin{enumerate}[label=$(\roman*)$]\itemsep=0pt
\item if $d+m\notin 2\mathbb{N}$, then $\big(\widetilde{\mathrm{Ric}_\phi^m},\widetilde{F_\phi^m}\big)=O_{IJ}(\rho^{\infty})$;

\item if $d+m\in 2\mathbb{N}$, then $\big(\widetilde{\mathrm{Ric}_\phi^m},\widetilde{F_\phi^m}\big)=O^{2,+}_{IJ}\big(\rho^{\frac{d+m}{2}-1}\big).$
 \end{enumerate}
\end{Definition}

A key step in the proof of Theorem~\ref{ambient-statement} is the construction and classification of straight and normal weighted ambient metrics \cite{FeffermanGraham2012}.

\begin{Theorem}
 \label{weighted_ambient_metrics_exist}
 Let $(M^d,g,f,m,\mu)$, $d\geq 3$ and $m<\infty$, be a smooth metric measure space. Then there exists a straight and normal weighted ambient space $(\cmG,\cg,\cf,m,\mu)$ for it. Moreover, if $(\cmG_j,\cg_j,\cf_j,m,\mu)$, $j\in\{1,2\}$, are two such weighted ambient spaces, then
\begin{enumerate}[label=$(\roman*)$]\itemsep=0pt
 \item if $d+m\notin 2\mathbb{N}$, then $\big(\cg_1-\cg_2,-2(\cf_1-\cf_2)\big)=O(\rho^\infty)$; and
 \item if $d+m\in 2\mathbb{N}$, then $\big(\cg_1-\cg_2,-2(\cf_1-\cf_2)\big)=O_{IJ}^{1,+}\big(\rho^{\frac{d+m}{2}-1}\big)$.
 \end{enumerate}
\end{Theorem}

We prove Theorem~\ref{weighted_ambient_metrics_exist} in Section~\ref{formal-theory}.

Note that if $(\cg_1-\cg_2,\cf_1-\cf_2)=O_{IJ}\big(\rho^{\frac{d+m}{2}-1}\big)$, then
\begin{gather*}
{\partial^{\frac{d+m}{2}}_\rho}\big|_{\rho=0}\big(\widetilde{f}_1^m\mathrm{dvol}_{\widetilde{g}_1}- \widetilde{f}_2^m\mathrm{dvol}_{\widetilde{g}_2}\big)
\\ \qquad
{}= \biggl(mf^{-1}{\partial_\rho^{\frac{d+m}{2}}}\big|_{\rho=0}\big(\widetilde{f}_1-\widetilde{f}_2\big)
+\frac{1}{2}g^{ij}{\partial_\rho^{\frac{d+m}{2}}}\big|_{\rho=0} (\widetilde{g}_1-\widetilde{g}_2)_{ij}\biggr)f^m\mathrm{dvol}_{g}.
\end{gather*}
Thus, a geometric interpretation of the improved orders of vanishing in the case that $d+m\in 2\mathbb{N}$ is that the weighted volume element $\widetilde{f}^m\operatorname{dvol}_{\widetilde{g}}$ vanishes to one higher order in $\rho$ than the weighted ambient space $(\widetilde{g},\widetilde{f})$.

The following definition generalizes the notion of uniqueness from Theorem~\ref{weighted_ambient_metrics_exist} to general pre-ambient spaces.

\begin{Definition}
{\sloppy
We say that two weighted pre-ambient spaces $(\widetilde{\mathcal{G}}_{1}, \widetilde{g}_{1},\widetilde{f}_1,m,\mu)$ and $(\widetilde{\mathcal{G}}_{2}, \widetilde{g}_{2}, \widetilde{f}_2,\allowbreak m,\mu)$ for $(M^d,[g,f],m,\mu)$ are \emph{ambient-equivalent} if there exist open sets $\mathcal{U}_{1} \subset \widetilde{\mathcal{G}}_{1}$ and $\mathcal{U}_{2} \subset \widetilde{\mathcal{G}}_{2}$, and a diffeomorphism $\phi\colon \mathcal{U}_{1} \rightarrow \mathcal{U}_{2}$ such that}
\begin{enumerate}[label=$(\roman*)$]\itemsep=0pt
\item $\mathcal{U}_{1}$ and $\mathcal{U}_{2}$ both contain $\mathcal{G} \times\{0\}$;

\item $\mathcal{U}_{1}$ and $\mathcal{U}_{2}$ are dilation-invariant and $\phi$ commutes with dilations;

\item the restriction of $\phi$ to $\mathcal{G} \times\{0\}$ is the identity map; and

\item \begin{itemize}
\item[$(a)$] if $d+m\notin 2\mathbb{N}$, then $\big(\widetilde{g}_{1}-\phi^{*} \widetilde{g}_{2}, -2(\widetilde{f}_1-\phi^*\widetilde{f}_2)\big)=O(\rho^\infty)$;
\item[$(b)$] if $d+m\in 2\mathbb{N}$, then $\big(\widetilde{g}_{1}-\phi^{*} \widetilde{g}_{2}, -2(\widetilde{f}_1-\phi^*\widetilde{f}_2)\big)=O_{I J}^{1,+}\big(\rho^{\frac{d+m}{2}-1}\big)$.
\end{itemize}
\end{enumerate}
\end{Definition}

Using Theorem~\ref{normal-to-straight-and-normal} below, we see that every weighted pre-ambient space is ambient-equivalent to a straight and normal weighted pre-ambient space.

{\sloppy\begin{Theorem}\label{pullback-normal}
Let $(\cmG,\cg,\cf,m,\mu)$ be a weighted ambient space for the conformal class $(M^d,[g,f],m,\mu)$, $d \geq 3$ and $m<\infty$, of smooth metric measure spaces. Then it is ambient-equivalent to a weighted ambient space in straight and normal form relative to $(g,f)$.
\end{Theorem}}

\begin{proof}
Following the proof in \cite[Proposition 2.8]{FeffermanGraham2012}, we observe that any weighted pre-ambient space is ambient-equivalent to a normal weighted pre-ambient space. Moreover, it follows from Theorem~\ref{normal-to-straight-and-normal} below that a normal weighted ambient space is ambient-equivalent to a straight and normal ambient space through the identity map.
\end{proof}

The uniqueness of Theorem~\ref{weighted_ambient_metrics_exist} implies the uniqueness of weighted ambient spaces up to ambient-equivalence (cf.\ \cite[Theorem~2.9]{FeffermanGraham2012}).

\begin{Theorem}
Any two weighted ambient spaces for $(M,[g,f],m,\mu)$, $d\geq 3$ and $m<\infty$, are ambient-equivalent.
\end{Theorem}

{\sloppy\begin{proof}
We pick a representative $(g,f)$ and invoke Theorem~\ref{weighted_ambient_metrics_exist}. Then there exists a straight and normal weighted ambient space $(\widetilde{\mathcal{G}}_1,\widetilde{g}_1,\widetilde{f}_1,m,\mu)$ for $(M,g,f,m,\mu)$. Now let $(\widetilde{\mathcal{G}}_2,\widetilde{g}_2,\widetilde{f}_2,m,\mu)$ be a~weighted ambient space for $(M,[g,f],m,\mu)$. Applying Theorem~\ref{pullback-normal}, we find that $(\widetilde{\mathcal{G}}_2, \widetilde{g}_2,\widetilde{f}_2,m,\mu)$ is ambient-equivalent to a weighted ambient space in straight and normal form relative to $(g,f)$. By Theorem~\ref{weighted_ambient_metrics_exist}, this space is ambient-equivalent to $(\widetilde{\mathcal{G}}_1,\widetilde{g}_1,\widetilde{f}_1,m,\mu)$.
\end{proof}}

\section{Formal theory}\label{formal-theory}

We prove Theorem~\ref{weighted_ambient_metrics_exist} by iteratively constructing a power series solution to $\widetilde{R}(\widetilde{g}),\widetilde{F_\phi^m}=O(\rho^j)$. This process is only obstructed when $d+m\in 2\mathbb{N}$; in this case the obstruction is at order $O\big(\rho^{\frac{d+m}{2}-1}\big)$.

We first give a necessary condition for a weighted ambient space to be normal.

\begin{Lemma}\label{ambient-metric-base-case}
Let $(\widetilde{\mathcal{G}},\widetilde{g},\widetilde{f},m,\mu)$ be a normal weighted ambient space. Then \begin{equation}\label{standard-normal-form}\widetilde{g}=a\,\mathrm{d}t^2+2b_i\, {\rm d}x^i \,\mathrm{d}t+ 2t\,{\rm d}\rho\, \mathrm{d}t+ t^2 (g_{ij})_\rho,\end{equation} with $a=2\rho+O(\rho^2)$ and $b_i=O(\rho^2)$.
\end{Lemma}

\begin{proof}
Lemma~\ref{straight-metric-lemma} implies that $\widetilde{g}$ has the form of equation~\eqref{standard-normal-form}. Since $\widetilde{g}$ is normal, $a|_{\rho=0}=0$ and $b_i|_{\rho=0}=0$. Denote $\widetilde{R}_{IJ}:=\widetilde{(\ric)}_{IJ}$.
By definition,
\begin{gather*}
\begin{split}
&\widetilde{R}_{IJ}=\frac{1}{2}\widetilde{g}^{KL}\big(\partial^2_{IL}\widetilde{g}_{JK} +\partial^2_{JK}\widetilde{g}_{IL}-\partial^2_{KL}\widetilde{g}_{IJ}-\partial^2_{IJ}\widetilde{g}_{KL}\big)
\\ &\hphantom{\widetilde{R}_{IJ}=}
{}+\widetilde{g}^{KL}\widetilde{g}^{PQ}\big(\widetilde{\Gamma}_{ILP}\widetilde{\Gamma}_{JKQ} -\widetilde{\Gamma}_{IJP}\widetilde{\Gamma}_{KLQ}\big)
-\frac{m}{\widetilde{f}}\big(\partial^2_{IJ}\widetilde{f} -\widetilde{\Gamma}^K_{IJ}\partial_K\widetilde{f}\big).
\end{split}
\end{gather*}
Evaluating at $\rho=0$ yields
\begin{gather*}
\widetilde{R}_{00}= \frac{d+m}{2t^2}(2-\partial_{\rho}a),
\qquad
\widetilde{R}_{0i}=\frac{1}{2t}\partial^2_{i\rho}a -\frac{d+m}{2t}\partial_{\rho}b_i.
\end{gather*}
We conclude that $a=2\rho+O(\rho^2)$ and $b_i=O(\rho^2)$.
\end{proof}

We now take a brief digression to discuss straight weighted ambient spaces. First, notice that a simple choice of certain components of $\widetilde{g}$ makes $\widetilde{R}_{0I}$ vanish.

\begin{Lemma}
\label{straight-zero-theorem}
Let $(\widetilde{\mathcal{G}},\widetilde{g},\widetilde{f},m,\mu)$ be a weighted pre-ambient space. If $\widetilde{g}$ has the form
\begin{equation}
\label{straightmetric}
 \widetilde{g}_{IJ}=\begin{pmatrix} 2\rho&0&t\\0& t^2 g_{\rho}&0\\t&0&0 \end{pmatrix}\!,
\end{equation}
where $g_\rho$ is a one-parameter family of metrics on $M$, then $\widetilde{R}_{0I}=0$.
\end{Lemma}

\begin{proof}This follows by direct computation~(cf.\ \cite[Lemma 3.2]{FeffermanGraham2012}).
\end{proof}

The relevance of Lemma~\ref{straight-zero-theorem} stems from the following equivalent characterizations of straight normal pre-ambient metrics.
\begin{Proposition}
\label{straightproperties}
Let $(\widetilde{\mathcal{G}},\widetilde{g},\widetilde{f},m,\mu)$, $d\geq 3$ and $m<\infty$, be in normal form relative to $(g,f)$. Then the following conditions are equivalent:
\begin{enumerate}[label=$(\roman*)$]\itemsep=0pt
\item $\widetilde{g}_{00}=2\rho$ and $\widetilde{g}_{0i}=0$;

\item for each $p\in \widetilde{\mathcal{G}}$, the dilation orbit $s\mapsto \delta_s p$ is a geodesic for $\widetilde{g}$;

\item $\widetilde{g}(2T,\cdot)={\rm d}(\widetilde{g}(T,T))$;

\item the infinitesimal dilation field $T$ satisfies $\widetilde{\nabla}T=\operatorname{Id}$.
\end{enumerate}
\end{Proposition}

\begin{proof}
The proof is identical to that of \cite[Proposition 3.4]{FeffermanGraham2012}.
\end{proof}

The following result implies that, up to ambient-equivalence, we may restrict our attention to metrics of the form of equation~\eqref{straightmetric}.
\begin{Theorem}
\label{normal-to-straight-and-normal}
A normal weighted ambient space is ambient-equivalent to a straight and normal ambient space.
\end{Theorem}
\begin{proof}
Let $(\widetilde{\mathcal{G}},\widetilde{g},\widetilde{f},m,\mu)$ be a normal weighted ambient space. By Lemma~\ref{ambient-metric-base-case}, we may write
\[
\widetilde{g}=\begin{pmatrix} 2\rho&0&t\\0&t^2 g_\rho& 0\\ t&0&0\end{pmatrix}+\rho^n \begin{pmatrix} a& tb_i& 0\\ tb_j & 0 & 0\\ 0&0&0\end{pmatrix}
\]
for some $n\geq 2$ and some $a=a(x,\rho)$ and $b_i=b_i(x,\rho)$. Direct computation using Lemma~\ref{straight-zero-theorem} yields
\begin{gather*}
t^{2} \widetilde{R}_{00}=n\bigg(n-1-\frac{d+m}{2}\bigg) \rho^{n-1} a+O(\rho^{n}),
\\
t \widetilde{R}_{0 i}=n\bigg(n-1-\frac{d+m}{2}\bigg) \rho^{n-1} b_i+\frac{n}{2} \rho^{n-1} \partial_{i} a+O(\rho^{n}),
\end{gather*}
(cf.\ \cite[equation~(3.11)]{FeffermanGraham2012}).
We conclude that if $d+m\notin 2\mathbb{N}$, then $a,b_i=O(\rho^\infty)$; and if $d+m\in 2\mathbb{N}$, then $n\geq \frac{d+m}{2}$. Therefore $\widetilde{g}$ is ambient-equivalent to a metric of the form of equation~\eqref{straightmetric}. The~conclusion follows from Proposition~\ref{straightproperties}.
\end{proof}

We now iteratively construct a straight and normal weighted ambient space.
Let $n\in\mathbb{N}$ be such that there is a metric
\[
\widetilde{g}_{IJ}^{(n-1)} = \begin{pmatrix}
2\rho & 0 & t \\
0 & t^2 g_{\rho} & 0 \\
t & 0 & 0
\end{pmatrix}
\]
and a function $\widetilde{f}^{(n-1)}=tf_\rho$ such that
\begin{gather*}
\widetilde{R}^{(n-1)}=O\big(\rho^{n-1}\big),
\qquad
\widetilde{F_{\phi}^m}^{(n-1)}=O\big(\rho^{n-1}\big).
\end{gather*}
Note that the existence of $\widetilde{g}_{IJ}^{(1)}$ and $\widetilde{f}^{(1)}$ trivially holds. Let $\cg^{(n)}_{IJ}=\cg^{(n-1)}_{IJ}+\Phi_{IJ}$ and $\widetilde{f}^{(n)}= \widetilde{f}^{(n-1)}+\rho^{n}t\upsilon$, where
\begin{align*}
\Phi_{IJ}&= \rho^n \begin{pmatrix}
0 & 0 & 0 \\
0 & t^2\psi_{ij} & 0 \\
0 & 0 & 0
\end{pmatrix}
\end{align*}
and $\psi_{ij}$, $\upsilon$ depend only on $x$ and $\rho$.
We seek $\psi_{IJ}$ and $\upsilon$ such that
 \begin{gather*}
\widetilde{R}^{(n)}= O(\rho^{n}),
\qquad
\widetilde{F_{\phi}^m}^{(n)} = O(\rho^{n}).
\end{gather*}
Note that the inverse of $\cg^{(n)}$, calculated modulo $O(\rho^n)$, is
\[
\cg^{IJ} = \begin{pmatrix} 0 & 0 & t^{-1}\\ 0 & t^{-2}g^{ij}& 0\\ t^{-1} & 0 & -{2\rho}t^{-2}\end{pmatrix}\!.
\]

The Christoffel symbols for $\widetilde{g}^{(n)}$ modulo $O(\rho^n)$ are (cf.~\cite[equation~(3.16)]{FeffermanGraham2012}).
\begin{gather*}
\widetilde{\Gamma}_{I J }^0=\begin{pmatrix}
0 & 0 & 0 \\
0 & -\frac{1}{2} t\big(\partial_{\rho} g_{i j}+n \rho^{n-1} \psi_{i j}\big) & 0 \\
0 & 0 & 0
\end{pmatrix}\!,
\\
\widetilde{\Gamma}_{I J}^\ell=\begin{pmatrix}
0 & t^{-1} \delta_{j}^{\ell} & 0 \\
t^{-1} \delta_{i}^{\ell} & \Gamma_{i j}^\ell & \frac{1}{2} g^{\ell k} \partial_{\rho} g_{i k}+\frac{n}{2} \rho^{n-1} \psi_{i}^{\ell}\\
0 & \frac{1}{2} g^{\ell k} \partial_{\rho} g_{j k}+\frac{n}{2} \rho^{n-1} \psi_{j}^{\ell} & 0
\end{pmatrix}\!,
\\
\widetilde{\Gamma}_{I J}^\infty=\begin{pmatrix}
0 & 0 & t^{-1} \\
0 & \rho \partial_{\rho} g_{i j}-g_{i j} & 0 \\
t^{-1} & 0 & 0
\end{pmatrix}\!.
\end{gather*}
It follows that~(cf.~\cite[equation~(3.11)]{FeffermanGraham2012})
\begin{subequations}
\begin{gather}
\label{ricij}
\widetilde{R}_{ij}^{(n)}=\widetilde{R}_{ij}^{(n-1)}+n\rho^{n-1} \bigg[\bigg(n-\frac{d+m}{2}\bigg)\psi_{ij}-\frac{1}{2}g^{kl}\psi_{kl}g_{ij}-\frac{m}{f}\upsilon g_{ij}\bigg]+O\big(\rho^{n}\big),
\\
\label{rphi}
\widetilde{F_{\phi}^m}^{(n)}= \widetilde{F_{\phi}^m}^{(n-1)}
+ n\rho^{n-1}f \bigg[\frac{1}{2}fg^{kl}\psi_{kl}+\upsilon (d+2m-2n) \bigg]+O\big(\rho^{n}\big),
\\
\widetilde{R}_{\infty\infty}^{(n)}=\widetilde{R}_{\infty\infty}^{(n-1)}-n(n-1)\rho^{n-2} \bigg[\frac{1}{2}g^{kl}\psi_{kl}+mf^{-1}\upsilon\bigg]+O\big(\rho^{n-1}\big)\label{ricinftyinfty}.
\end{gather}
\end{subequations}

\subsection[Solving widetilde {R}\_\{ij\}, widetilde {F\_\{phi\}\textasciicircum{}m}=O(rho\textasciicircum{}n)]{Solving $\boldsymbol{\widetilde{R}_{ij},\widetilde{F_\phi^m}=O(\rho^n)}$}
\label{section-4.1}

The following theorem specifies the result of recursively determining $(\widetilde{g}^{(n)},\widetilde{f}^{(n)})$ by the requirements $\widetilde{R}_{ij}^{(n)},(\widetilde{F}_\phi^m)^{(n)}=O(\rho^n)$.
\begin{Theorem}
\label{various-possibilities}
Let $(M^d,g,f,m,\mu)$, $d\geq 3$ and $m<\infty$, be a smooth metric measure space. Set $\widetilde{\mathcal{G}}:=\mathbb{R}_+\times M\times (-\epsilon,\epsilon)$.
\begin{enumerate}[label=$(\roman*)$]\itemsep=0pt
\item If $d+m\notin\mathbb{N}$, then there is a straight and normal metric measure structure $(\widetilde{g},\widetilde{f})$ on $\widetilde{\mathcal{G}}$, unique modulo $O(\rho^\infty)$, such that $\widetilde{R}_{ij},\widetilde{F_\phi^m}=O(\rho^\infty)$.

\item If $d+m\in 2\mathbb{N}$, then
\begin{enumerate}\itemsep=0pt
\item[$(a)$] there is a straight and normal metric measure structure $(\widetilde{g},\widetilde{f})$ on $\widetilde{\mathcal{G}}$, unique modulo $O^+\big(\rho^{\frac{d+m}{2}}\big)$, such that $(\widetilde{R}_{ij},\widetilde{F_\phi^m})=O_{ij}^{2,+}\big(\rho^{\frac{d+m}{2}-1}\big)$;
\item[$(b)$] if
\[
\partial_\rho^{(d+m)/2-1}\big|_{\rho=0}\widetilde{R}_{ij}=0,
\]
then there is a straight and normal metric measure structure $(\widetilde{g},\widetilde{f})$ on $\widetilde{\mathcal{G}}$ such that $\widetilde{R}_{ij},\widetilde{F_\phi^m}=O(\rho^{d+m-1})$;

\item[$(c)$] if
\begin{equation}\label{oddconsistency}
\frac{m}{f^2} \widetilde{F_\phi^m}^{(d+m-1)}-g^{ij}\widetilde{R}_{ij}^{(d+m-1)}=O\big(\rho^{d+m}\big),
\end{equation}
then there is a straight and normal metric measure structure $(\widetilde{g},\widetilde{f})$ on $\widetilde{\mathcal{G}}$ such that $\widetilde{R}_{ij},\widetilde{F_\phi^m}=O(\rho^\infty)$.
\end{enumerate}

\item If $d+m\in \mathbb{N}\setminus 2\mathbb{N}$, then
\begin{enumerate}\itemsep=0pt
\item[$(a)$] there is a straight and normal metric measure structure $(\widetilde{g},\widetilde{f})$ on $\widetilde{\mathcal{G}}$, unique modulo $O(\rho^{d+m})$, such that $\widetilde{R}_{ij},\widetilde{F_\phi^m}=O(\rho^{d+m-1})$;
\item[$(b)$] if equation~\eqref{oddconsistency} holds, then there is a straight and normal metric measure structure $(\widetilde{g},\widetilde{f})$ on $\widetilde{\mathcal{G}}$ such that $\widetilde{R}_{ij},\widetilde{F_\phi^m}=O(\rho^\infty)$.
\end{enumerate}
\end{enumerate}
\end{Theorem}

\begin{Remark}
We say that $(\widetilde{g},\widetilde{f})$ is unique modulo $O^+(\rho^k)$ if it is unique modulo $O(\rho^{k})$ and $\frac{1}{2}g^{ij}\widetilde{g}_{ij}+\frac{m}{f}\widetilde{f}$ is unique modulo $O(\rho^{k+1})$.
\end{Remark}

\begin{proof}
On studying equations~\eqref{ricij} and \eqref{rphi}, we observe that if $n\neq \frac{d+m}{2}$, then there is a~unique choice of the trace-free part $\psi_{ij}-\frac{1}{d}g^{kl}\psi_{kl}g_{ij}$ of $\psi_{ij}$ which makes the trace-free part of $\widetilde{R}^{(n)}_{ij}$ vanish modulo $O(\rho^n)$. However, if $n=\frac{d+m}{2}$, then we may not be able to make the trace-free part of $\widetilde{R}_{ij}^{(n)}$ vanish modulo $O(\rho^n)$. Additionally, equations~\eqref{ricij} and \eqref{rphi} imply that
\begin{gather}
g^{ij}\widetilde{R}_{ij}^{(n)}=g^{ij}\widetilde{R}_{ij}^{(n-1)}
+n\rho^{n-1}\frac{(2n-2d-m)}{2}g^{kl}\psi_{kl}-\frac{dmn}{f}\upsilon \rho^{n-1}+O\big(\rho^{n}\big),\nonumber
\\
\widetilde{F_{\phi}^m}^{(n)}= \widetilde{F_{\phi}^m}^{(n-1)}
+ n\rho^{n-1}f \bigg[\frac{1}{2}fg^{kl}\psi_{kl}+\upsilon (d+2m-2n) \bigg]+O\big(\rho^{n}\big).
\label{modified-equations}
\end{gather}
The determinant of the coefficient matrix of $g^{kl}\psi_{kl}$ and $\upsilon$ is
\begin{equation*}(2n-d-m)(n-d-m).\end{equation*}
Therefore if $n\notin\{\frac{d+m}{2},d+m\}$, then there are unique $g^{kl}\psi_{kl}$ and $\upsilon$ such that $g^{ij}\widetilde{R}_{ij}^{(n)},(\widetilde{F_\phi^m})^{(n)}=O(\rho^n)$. However, if $n\in\{\frac{d+m}{2},d+m\}$, then we may not be able to simultaneously solve $g^{ij}\widetilde{R}_{ij}^{(n)}=O(\rho^n)$ and $\widetilde{F_\phi^m}=O(\rho^n)$.

Suppose first that $d+m\notin \mathbb{N}$. Combining the above discussion with Borel's lemma yields a~straight and normal metric measure structure $(\widetilde{g},\widetilde{f})$ such that $\widetilde{R}_{ij},\widetilde{F_\phi^m}=O(\rho^\infty)$.

Suppose next that $d+m\in 2\mathbb{N}$. Set $n=\frac{d+m}{2}$. Equations~\eqref{modified-equations} imply that
\begin{align}
\frac{m}{f^2}\widetilde{F_{\phi}^m}^{(\frac{d+m}{2})}- g^{ij}\widetilde{R}_{ij}^{(\frac{d+m}{2})} ={}&\frac{m}{f^2}\widetilde{F_{\phi}^m}^{(\frac{d+m}{2}-1)}-g^{ij}\widetilde{R}_{ij}^{(\frac{d+m}{2}-1)} \nonumber
\\
&+\frac{(d+m)^2}{2}\rho^{\frac{d+m}{2}-1} \bigg(\frac{1}{2}g^{kl}\psi_{kl}+\frac{m}{f}\upsilon\bigg)+O\big(\rho^{\frac{d+m}{2}}\big).
\label{even-difference}
\end{align}
In particular, there is a unique choice of $\frac{1}{2}g^{kl}\psi_{kl}+\frac{m}{f}\upsilon$ such that the left hand side of equation~\eqref{even-difference} is $O(\rho^{\frac{d+m}{2}})$. Making this choice yields a unique straight and normal metric measure structure $(\widetilde{g},\widetilde{f})$ such that $(\widetilde{R}_{ij},\widetilde{F_\phi^m})=O_{ij}^{2,+}(\rho^{\frac{d+m}{2}-1})$.
Moreover, if the weighted obstruction is zero, then $\widetilde{R}_{ij}^{(n)}=O(\rho^{\frac{d+m}{2}})$. Hence, we may iteratively solve equations~\eqref{ricij} and \eqref{rphi} to obtain a straight and normal metric measure structure $(\widetilde{g},\widetilde{f})$ such that $\widetilde{R}_{ij},\widetilde{F_\phi^m}=O(\rho^{d+m-1})$. The possible obstruction is the same as that discussed in the next case.

Suppose finally that $d+m\in \mathbb{N}\setminus 2\mathbb{N}$. Set $n=d+m$. From the discussion of the first paragraph, we may choose the trace-free part of $\psi_{ij}$ such that $\widetilde{R}_{ij}^{(d+m)}$ is pure trace. Additionally, by equations~\eqref{modified-equations},
\begin{gather*}
g^{ij}\widetilde{R}_{ij}^{(d+m)}=g^{ij}\widetilde{R}_{ij}^{(d+m-1)}+
m(d+m) \rho^{d+m-1}\bigg(\frac{1}{2}g^{kl}\psi_{kl}-\frac{d}{f}\upsilon\bigg)+O\big(\rho^{d+m}\big),
\\
\widetilde{F_{\phi}^m}^{(d+m)}= \widetilde{F_{\phi}^m}^{(d+m-1)}
+f^2(d+m)\rho^{d+m-1} \bigg(\frac{1}{2}g^{kl}\psi_{kl}-\frac{d}{f}\upsilon\bigg)+O\big(\rho^{d+m}\big).
\end{gather*}
Therefore we may choose $g^{ij}\psi_{ij}$ and $\upsilon$ such that $g^{ij}\widetilde{R}_{ij}^{(d+m)},\widetilde{F_{\phi}^m}^{(d+m)}=O(\rho^{d+m})$ if and only if
equation~\eqref{oddconsistency} holds. If equation~\eqref{oddconsistency} holds, then we may continue iteratively improving $(\widetilde{g},\widetilde{f})$ as in the first paragraph. Hence, by Borel's lemma, there is a straight and normal metric measure structure $(\widetilde{g},\widetilde{f})$ such that $\widetilde{R}_{ij},\widetilde{F_\phi^m}=O(\rho^\infty)$.
\end{proof}

\subsection[Solving widetilde {R}\_\{I infty\}=0]{Solving $\boldsymbol{\widetilde{R}_{I\infty}=0}$}
In this subsection, we show that if $(\widetilde{g},\widetilde{f})$ is as in Theorem~\ref{various-possibilities}, then the components $\widetilde{R}_{I\infty}$ can be made to vanish to the appropriate order and that equation~\eqref{oddconsistency} automatically holds. This proves Theorem~\ref{weighted_ambient_metrics_exist}.

\begin{proof}[Proof of Theorem~\ref{weighted_ambient_metrics_exist}]
Let $(\widetilde{g},\widetilde{f})$ be as in Theorem~\ref{various-possibilities}. Equation~\eqref{intro-weighted-bianchi} implies that
\begin{gather}
\widetilde{g}^{JK}\widetilde{\partial}_J\widetilde{R}_{KI} -\widetilde{g}^{JK}\widetilde{\Gamma}_{JK}^Q\widetilde{R}_{QI} -\widetilde{g}^{JK}\widetilde{\Gamma}_{JI}^Q\widetilde{R}_{KQ} -\widetilde{g}^{JK}\widetilde{R}_{IJ}\partial_K \widetilde{\phi}
\!-\frac{1}{2}\partial_I \widetilde{R_{\widetilde{\phi}}^m}-\frac{1}{f^2}\widetilde{F}_\phi^m \partial_I\widetilde{\phi}=0.\!
\label{bianchianalogue}
\end{gather}

Set
\begin{equation*}
n=\begin{cases}
\infty & \text{if}\quad d+m\notin\mathbb{N},\\
\frac{d+m}{2}-1&\text{if}\quad d+m\in 2\mathbb{N},\\
d+m-1& \text{if}\quad d+m\in \mathbb{N}\setminus 2\mathbb{N}.
\end{cases}
\end{equation*}
Taking $I=l$ and $I=\infty$ in equation~\eqref{bianchianalogue} and computing$\;\bmod \;O(\rho^{n})$ yields
\begin{subequations}
\label{bianchieqns}
\begin{gather}
\label{bianchi2}
[d+m-2-2\rho\partial_{\rho}] \widetilde{R}_{\infty l}-\rho g^{ks}g'_{sl}\widetilde{R}_{\infty k}+2\rho \widetilde{R}_{l\infty}\partial_\rho \phi+\rho \partial_l\widetilde{R}_{\infty\infty}=O\big(\rho^{n}\big),
\\
[d+m-2-\rho\partial_{\rho}] \widetilde{R}_{\infty\infty}
-\frac{1}{2}\partial_{\rho}\bigg(g^{ij}\widetilde{R}_{ij}-\frac{m}{f^2}\widetilde{F_\phi^m}\bigg)
+g^{ij}\big(\widetilde{\nabla}_\phi\big)_i\widetilde{R}_{j\infty}\nonumber
\\ \qquad
+\rho\big(2\partial_\rho\phi-g^{ij}g'_{ij}\big)\widetilde{R}_{\infty\infty}
=O(\rho^{n})\label{bianchi3}.
\end{gather}
\end{subequations}

First, suppose that $d+m\notin\mathbb{N}$. We know that $\widetilde{R}_{ij},\widetilde{F_\phi^m}=O(\rho^\infty)$. From equations~\eqref{bianchi2} and~\eqref{bianchi3} we conclude that $\widetilde{R}_{\infty l},\widetilde{R}_{\infty\infty}=O(\rho^\infty)$.

Second, suppose that $d+m\in \mathbb{N}\setminus 2\mathbb{N}$. Equation~\eqref{bianchieqns} implies that $\widetilde{R}_{\infty l},\widetilde{R}_{\infty\infty}=O(\rho^{n})$, and that $g^{ij}\widetilde{R}_{ij}-mf^{-2}\widetilde{F_\phi^m}=O(\rho^{n+1})$. This verifies equation~\eqref{oddconsistency}. Equation~\eqref{ricinftyinfty} gives us the unique value of $\frac{1}{2}g^{kl}\psi_{kl}+\frac{m}{f}\upsilon$ such that $\widetilde{R}_{\infty\infty}=O(\rho^{n})$. Hence, we can now uniquely determine the values of $g^{kl}\psi_{kl}$ and $\upsilon$, and solve $\widetilde{R}_{IJ},\widetilde{F_\phi^m}=O(\rho^\infty)$.

Third, suppose that $d+m\in 2\mathbb{N}$. Equation~\eqref{bianchi2} tells us that $\widetilde{R}_{\infty l}=O(\rho^n)$ and that $\widetilde{R}_{\infty\infty}=O(\rho^{n-1})$. Now recall that $g^{ij}\widetilde{R}_{ij}-mf^{-2}\widetilde{F_\phi^m}=O\big(\rho^{\frac{d+m}{2}}\big)$.
Hence, equation~\eqref{bianchi3} tells us that $\widetilde{R}_{\infty\infty}=O(\rho^n)$.

Finally, since the terms of $(\widetilde{g},\widetilde{f})$ are uniquely determined using the process described above, any two straight and normal weighted ambient structures will be ambient-equivalent to the orders stated in the theorem.
\end{proof}

\subsection{The weighted obstruction tensor}
We conclude this section by establishing some properties of the weighted obstruction tensor. To~that end, we first introduce notation for the corresponding term in the asymptotic expansion of~$\widetilde{F_\phi^m}$.
\begin{Definition}
We define
\[
\mathcal{F}:=c_{d+m}\partial_\rho^{\frac{d+m}{2}-1}\big|_{\rho=0}\widetilde{F_\phi^m},
\]
where $c_{d+m}$ is defined in equation~\eqref{obstruction-definition}.
\end{Definition}

There is a simple relationship between $\mathcal{F}$ and the trace and divergence of $\mathcal{O}_{ij}$ (cf.\ \cite{Case2014s}).

\begin{Theorem}
\label{obstruction-theorem}
Let $(M^d,g,f,m,\mu)$, $d\geq 3$ and $m<\infty$, be a smooth metric measure space with $d+m\in 2\mathbb{N}$. Then
\begin{enumerate}[label=$(\roman*)$]\itemsep=0pt
\item $\mathcal{O}_i^i=\frac{m}{f^2}\mathcal{F}$;
\item $\delta_\phi\mathcal{O}_{i}=\frac{1}{f^2}\mathcal{F}\partial_i\phi$; and
\item $\mathcal{O}_{ij}$ and $\mathcal{F}$ are local weighted conformal invariants of weight $2-d-m$.
\end{enumerate}
\end{Theorem}

\begin{proof}
We can use the recursive construction of $(\widetilde{g},\widetilde{f})$ to express $\mathcal{O}_{ij}$ and $\mathcal{F}$ as polynomials in~$g$, $g^{-1}$, $f$, $\operatorname{Rm}$ and $\nabla$.
Hence $\mathcal{O}_{ij}$ and $\mathcal{F}$ are local weighted invariants.

Since $g^{ij}\widetilde{R}_{ij}-mf^{-2}\widetilde{F_\phi^m}=O\big(\rho^{\frac{d+m}{2}}\big)$, we see that $\mathcal{O}^i_i=\frac{m}{f^2}\mathcal{F}$.

Writing equation~\eqref{bianchi2} modulo $O\big(\rho^{\frac{d+m}{2}}\big)$ and recalling that $\widetilde{R}_{\infty l}=O\big(\rho^{\frac{d+m}{2}-1}\big)$, we get
\[
\frac{1}{t^2}\widetilde{\nabla}^j \widetilde{R}_{jl}-\frac{1}{t^2}g^{ij}\widetilde{R}_{li}\partial_j \widetilde{\phi}-\frac{1}{f^2}\widetilde{F_\phi^m}\partial_l \widetilde{\phi}=O\big(\rho^{\frac{d+m}{2}}\big).
\]
Hence, $\delta_\phi\mathcal{O}_{i}=\frac{1}{f^2}\mathcal{F}\partial_i\phi$.

The conformal invariance of the weighted ambient space implies that $\mathcal{O}_{ij}$ is conformally invari\-ant of weight $2-d-m$. The conformal invariance of $\mathcal{F}$ follows from the identity ${\mathcal{O}_i^i=\frac{m}{f^2}\mathcal{F}}$.
\end{proof}

\subsection[The first few terms in the expansions of g\_\{rho\} and f\_\{rho\}]{The first few terms in the expansions of $\boldsymbol{g_\rho}$ and $\boldsymbol {f_\rho}$}
\label{sec:low-order-expansion}

Although it is difficult in general to compute the terms of $(\widetilde{g},\widetilde{f})$, we are able to compute the first few terms by hand.

Let $(\widetilde{g},\widetilde{f})$ be as in equation~\eqref{weighted-straight-normal-form}. Then
\begin{gather}
\widetilde{R}_{ij}=\rho g''_{ij}-\rho g^{kl}g'_{ik}g'_{jl}+\frac{1}{2}\rho g^{kl}g'_{kl}g'_{ij}+\rho\frac{m}{\widetilde{f}} g_{ij}'\widetilde{f}'
-\biggl(\frac{d+m}{2}-1\biggr)g'_{ij}-\frac{1}{2}g^{kl}g'_{kl}g_{ij}\nonumber
\\ \hphantom{\widetilde{R}_{ij}=}
{}-\frac{m}{\widetilde{f}}g_{ij}\widetilde{f}'+\big(\mathrm{Ric}_{\phi}^m\big)_{ij},\nonumber
\\
\widetilde{F_\phi^m}=-2\rho f f''\!-\rho ff'g^{ij}g'_{ij}\!-2(m-1)\rho(f')^2
 +\frac{1}{2}f^2g^{ij}g'_{ij}+(2m+d-2)ff'+F_\phi^m.\!
\label{secondorderricequation}
\end{gather}

Using these, we readily compute that
\begin{align*}
g_{ij}'(\cdot,0)=2 \big(P_\phi^m\big)_{ij},\qquad
f'(\cdot,0)=\frac{f}{m}Y_\phi^m.
\end{align*}
Moreover, differentiating equation~\eqref{secondorderricequation} once with respect to $\rho$ and evaluating at $\rho=0$ yields
\begin{gather*}
(d+m-4) g_{i j}^{\prime \prime}=-2\big(B_\phi^m\big)_{i j}+2(d+m-4) \big(P_\phi^m\big)_{i}^{k} \big(P_\phi^m\big)_{j k},
\\
(d+m-4)f''= \frac{f}{m}g^{ij}\big(B_\phi^m\big)_{ij}.
\end{gather*}
In particular, if $d+m=4$, then $\mathcal{O}_{ij}=(B_\phi^m)_{ij}$.

\section{Weighted Poincar\'e spaces}
\label{sec:poincare}

In this section, we construct a weighted Poincar\'e space for a smooth metric measure space.

Let $M_+$ be a smooth manifold with compact boundary $M$. Denote by $M^0_+$ the interior of $M_+$. A smooth metric measure space $(M^0_+,g_+,r_+,m,\mu)$, $m<\infty$, is \emph{conformally compact} if there is a defining function $r$ for $M$ such that
\begin{enumerate}[label=$(\roman*)$]\itemsep=0pt
 \item $(r^2g_+,rf_+)$ extends smoothly to $M_+$; and
 \item $({M,r^2 g_+}|_{M},{rf_+}|_{M},m,\mu)$ is a smooth metric measure space.
\end{enumerate} We call $(M,[{r^2 g_+}|_{M},{rf_+}|_{M}],m,\mu)$ the conformal infinity of $(M^0_+,g_+,r_+,m,\mu)$.

In the following, we identify $M_+$ with an open neighborhood of $M\times \{0\}$ in $M\times [0,\infty)$ and choose $r$ to be the standard coordinate on $[0,\infty)$.

\begin{Definition}
\label{poincare-definition}
A \emph{weighted Poincar\'e space} for $(M^d,[g,f],m,\mu)$, $m<\infty$, is a metric measure structure $(g_+,f_+)$ on $M\times [0,\epsilon)$ such that
\begin{enumerate}[label=$(\roman*)$]\itemsep=0pt

\item $g_+$ has signature $(d+1,0)$;

\item $(g_+,f_+)$ has $(M^d,[g,f],m,\mu)$ as conformal infinity; and

\item \begin{enumerate}\itemsep=0pt
\item[$(a)$] if $d+m\notin 2\mathbb{N}$, then
\begin{equation*}
\big(\mathrm{Ric}_{\phi}^m(g_+)+(d+m)g_+,F_\phi^m(g_+)
-(d+m)f_+^2\big)=O\big(r^\infty\big);
\end{equation*}

\item[$(b)$] if $d+m\in 2\mathbb{N}$, then
\begin{equation*}
\big(\mathrm{Ric}_{\phi}^m(g_+)+(d+m)g_+,F_\phi^m(g_+)
-(d+m)f_+^2\big)=O^{2,+}_{\alpha\beta}\big(r^{d+m-2}\big).
\end{equation*}
\end{enumerate}
\end{enumerate}
\end{Definition}

If $\psmms$ is conformally compact, then $|{\rm d}r/r|_{g_+}=|{\rm d}r|_{r^2 g_+}$ extends smoothly to the boundary. We call $\psmms$ \emph{asymptotically hyperbolic} if $|{\rm d}r/r|_{g_+}=1$ on $M$. In~this case, all the sectional curvatures of $g_+$ approach $-1$ at a boundary point~\cite{MazzeoMelrose1987}.

\begin{Definition}
A weighted Poincar\'e space $\psmms$, $m<\infty$, is in \emph{normal form relative to $(g,f)$} if $g_+=r^{-2}({\rm d}r^2+g_r)$ and $f_+=r^{-1}f_r$. Here $g_r$ is a one-parameter family of metrics on $M$ such that $g_0=g$, and $f_r$ is a one-parameter family of functions such that $f_0=f$.
\end{Definition}

\begin{Proposition}
Let $\psmms$ be asymptotically hyperbolic. Then there exists a~neigh\-bor\-hood $U$ of $M\times \{0\}$ in $M\times [0,\infty)$ on which there is a unique diffeomorphism $\phi$ from $U$ into~$M_+$ such that ${\phi}|_{M}$ is the identity map and $(\phi^* g_+,\phi^* f_+)$ is in normal form relative to $(g,f)$ on~$U$.
\end{Proposition}
\begin{proof}
This follows as in the conformal case \cite[Proposition 4.3]{FeffermanGraham2012}.
\end{proof}

\begin{Definition}
{\sloppy
An asymptotically hyperbolic smooth metric measure space, denoted as $\psmms$, is \emph{even} if $(r^2 g_+,rf_+)$ is the restriction to $M_+$ of a metric measure structure $(g_+',f_+')$ on an open set $V\subset M\times (-\infty,\infty)$ containing $M_+$ such that $V$ and $(g_+',f_+')$ are invariant under $r\mapsto -r$.

}

A diffeomorphism $\psi\colon M_+\mapsto M\times [0,\infty)$ satisfying ${\psi}"_{{M\times \{0\}}}=\mathrm{Id}$ is \emph{even} if $\psi$ is the restriction of a diffeomorphism on an open set $V$ as above which commutes with $r\mapsto -r$.
\end{Definition}

It is easily seen that if $\psi$ if an even diffeomorphism and $\psmms$ is an even asymptotically hyperbolic smooth metric measure, then $\psi^{*}g_+$ is also even.

\begin{Theorem}
\label{existence-even-poincare}
Let $(M^d,[g,f],m,\mu)$, $d\geq 3$ and $m<\infty$, be a conformal class of smooth metric measure spaces. Then there exists an even weighted Poincar\'e space corresponding to it. Moreover, if $(g_+^1,f_+^1)$ and $(g_+^2,f_+^2)$ are two even weighted Poincar\'e structures defined on $(M_+^1)^{\circ}$ and $(M_+^2)^{\circ}$ respectively, then there exist open sets $U_1\subset M_+^1$ and $U_2\subset M_+^2$ containing $M\times \{0\}$, and an even diffeomorphism $\phi\colon U_1\mapsto U_2$, such that ${\phi}|_{{M\times \{0\}}}$ is the identity map and
\begin{enumerate}[label=$(\roman*)$]\itemsep=0pt

\item if $d+m\notin 2\mathbb{N}$, then $(g_+^1-\phi^*{g_+^2},-2(f_+^1-\phi^*f_+^2))=O(\rho^\infty)$;

\item if $d+m\in 2\mathbb{N}$, then $(g_+^1-\phi^*{g_+^2},-2(f_+^1-\phi^*f_+^2))=O_{\alpha\beta}^+(r^{d+m-2})$.
\end{enumerate}
\end{Theorem}

Theorem~\ref{existence-even-poincare} will be proved at the end of this section.

Let $(\widetilde{\mG},\widetilde{g},\widetilde{f},m,\mu)$ be a straight weighted pre-ambient space for $(M,g,f,m,\mu)$. Then $\widetilde{g}(T,T)$ vanishes to first order on $\mG\times \{0\}\subset \widetilde{\mG}$. Therefore, shrinking $\widetilde{\mG}$ if necessary, the hypersurface $H:=\widetilde{\mG}\cap \{\widetilde{g}(T,T)=-1\}$ lies on one side of $\mG\times \{0\}$. Also, since $\widetilde{g}(T,T)$ is homogeneous of degree $2$, each $\delta_s$-orbit on this side of $\mathcal{G}\times \{0\}$ intersects $H$ exactly once. We extend the projection $\pi\colon \mG\mapsto M$ to $\widetilde{\mG}\subset \mG\times \mathbb{R}\mapsto M\times \mathbb{R}$ by projecting only the first factor. Define $\chi\colon M\times \mathbb{R}\mapsto M\times [0,\infty)$ by $(x,\rho)\mapsto (x,\sqrt{2|\rho|})$. Then there is an open set $M_+\subset M\times [0,\infty)$ containing $M\times \{0\}$ such that $\chi\circ \pi|_{H}\colon H\mapsto M_+^{\circ}$ is a diffeomorphism.

\begin{Proposition}
\label{straight-poincare-existence}
If $(\widetilde{\mathcal{G}},\widetilde{g},\widetilde{f},m,\mu)$ is a straight and normal weighted pre-ambient space of the form of equation~\eqref{weighted-straight-normal-form}, and if $H$ and $M_+^{\circ}$ are as above, then
\begin{gather*}
g_+:= \big((\chi\circ \pi|_{H})^{-1}\big)^* \widetilde{g},\qquad
f_+:= \big((\chi\circ \pi|_{H})^{-1}\big)^* \widetilde{f},
\end{gather*}
defines an even, asymptotically hyperbolic, normal smooth metric measure space $(M_+,g_+,f_+,\allowbreak m,\mu)$ with conformal infinity $(M^d,g,f,m,\mu)$.
\end{Proposition}

\begin{proof}
We have $H=\{2\rho t^2=-1\}$. Introduce new variables $r>0$ and $s>0$ by $-2\rho=r^2$ and $s=rt$. Then
\begin{gather*}
\widetilde{g}= \frac{s^2}{r^2}\big(g_{-\frac{1}{2}r^2}+{\rm d}r^2\big)-{\rm d}s^2,
\qquad
\widetilde{f}= \frac{s}{r}f_{-\frac{1}{2}r^2},
\end{gather*}
Note that $H=s^{-1}(\{1\})$. Thus the restriction of $\widetilde{g}$ to $TH$ is $r^{-2}(g_{-\frac{1}{2}r^2}+{\rm d}r^2)$. Similarly, the restriction of $\widetilde{f}$ to $H$ is $r^{-1}f_{-\frac{1}{2}r^2}$. By the definition of $\chi$, we see that $r$ is the coordinate in the second factor of $M\times [0,\infty)$. Therefore $(g_+,f_+)$ is an even asymptotically hyperbolic smooth metric measure space with conformal infinity $(M^d,[g,f],m,\mu)$ in normal form relative to $(g,f)$.
\end{proof}

We now show that weighted Poincar\'e spaces and weighted ambient spaces are closely related.

\begin{Proposition}
\label{ricci-poincare-condition}
Let $(M_+^{d+1},g_+,f_+,m,\mu)$ be a smooth metric measure space. Consider $(M_+\times \mathbb{R}_+,\widetilde{g},\widetilde{f},m,\mu)$, where $\widetilde{g}=s^2 g_+-{\rm d}s^2$ and $\widetilde{f}=sf_+$. Then
\begin{gather*}
\mathrm{Ric}_{\phi}^m(\widetilde{g})=\mathrm{Ric}(g_+)+(d+m)g_+,
\qquad
\widetilde{F_\phi^m} =F_\phi^m(g_+)-(d+m)f_+^2.
\end{gather*}
\end{Proposition}

\begin{proof}
The Christoffel symbols of the metric $\widetilde{g}=-{\rm d}s^2+s^2g_+$ are
\begin{equation}
\label{poincare-christoffel-symbols}
\widetilde{\Gamma}_{IJ}^0=\begin{pmatrix}0&0\\0& s(g_+)_{ij} \end{pmatrix}\!,\qquad
\widetilde{\Gamma}_{IJ}^k=\begin{pmatrix}0&\frac{1}{s}\delta_i^k\vspace{1mm}\\\frac{1}{s}\delta_i^k& (\Gamma_{g_+})_{ij}^k \end{pmatrix}\!.
\end{equation}
Using the formula
\begin{equation*}
\mathrm{Ric}_{i j}=\partial_L \Gamma_{i j}^{L}-\partial_i \Gamma_{L j}^{L}+\Gamma_{i j}^{P} \Gamma_{L P}^{L}-\Gamma_{L j}^{P} \Gamma_{i P}^{L}
\end{equation*}
and equation~\eqref{poincare-christoffel-symbols}, we readily compute that $\Ric_{\phi}^m(\widetilde{g})=\Ric(g_+)+(d+m)g_+$.
Using the formula
\begin{equation*}
 \widetilde{F_\phi^m}=sf_+\biggl(-\frac{d+1}{s}f_++\frac{1}{s^2}\Delta_{g_+}sf_+\biggr) +(m-1)\big({-}f_+^2+|\nabla_{g_+}f_+|^2-\mu\big),
\end{equation*}
we readily compute that $\widetilde{F_\phi^m}=F_\phi^m(g_+)-(d+m)f_+^2.$ \qedhere
\end{proof}

\begin{proof}[Proof of Theorem~\ref{existence-even-poincare}]
Pick a representative $(M^d,g,f,m,\mu)$. Theorem~\ref{weighted_ambient_metrics_exist} implies that there is a weighted ambient metric in straight and normal form relative to $(g,f)$. By Proposition~\ref{straight-poincare-existence}, the metric measure structure $(g_+,f_+)$ is even and in normal form relative to $(g,f)$. Proposition~\ref{ricci-poincare-condition} implies that conditions (3) and (4) of Definition~\ref{poincare-definition} are satisfied, so that $(M_+^\circ,g_+,f_+,m,\mu)$ is a weighted Poincar\'e metric space for $(M^d,g,f,m,\mu)$. This proves the first part of the theorem. The second part follows similarly from the argument in Theorem~\ref{weighted_ambient_metrics_exist}.
\end{proof}

\section{Explicit examples of weighted ambient structures}
\label{sec:explicit-examples}

We conclude this article by proving Theorem~\ref{examples-theorem}. The following three subsections handle the three separate cases mentioned in the theorem.

\subsection{Quasi-Einstein space}

Suppose that $(M^d,g,f,m,\mu)$, $m<\infty,$ is quasi-Einstein. Then
\begin{equation*}
\mathrm{Ric}_{\phi}^m=\lambda g,
\qquad F_\phi^m=-\lambda f^2.
\end{equation*}
Set
\begin{equation*}
g_\rho:=(1+\lambda \rho)^2 g,\qquad
f_\rho:=(1+\lambda\rho)f.
\end{equation*}
Direct computation yields $\widetilde{\mathrm{Ric}_\phi^m}=0$ and $\widetilde{F_\phi^m}=0$. Thus $(\mathcal{G}\times (-\epsilon,\epsilon),\widetilde{g},\widetilde{f},m,\mu)$ is a weighted ambient metric for $(M^d,g,f,m,\mu)$. This verifies Theorem~\ref{examples-theorem}\ref{Theorem 1.2-1}.

\subsection{Weighted locally conformally flat space}

Recall from Section~\ref{sec:bg.tex} that the locally conformally flat condition is equivalent to $A_\phi^m=0$ and ${\rm d}P_\phi^m=0$.
Therefore $(P_\phi^m)_{ij;k}=(P_\phi^m)_{(ij;k)}$. We use this fact repeatedly in the rest of the section.

Let $(g_\rho,f_\rho)$ be as in the statement of Theorem~\ref{examples-theorem}. Note that
\[
(g_\rho)_{i j}=g_{i l} \big(U_\phi^m\big)_{k}^{l} \big(U_\phi^m\big)_{j}^{k}=\big(U_\phi^m\big)_{i k} \big(U_\phi^m\big)_{j}^{k},
\]
where
\[
\big(U_\phi^m\big)_{j}^{i}:=\delta_{j}^{i}+ \rho\big(P_\phi^m\big)_{j}^{i}.
\]
Set $V_\phi^m:=(U_\phi^m)^{-1}$. Then $(V_\phi^m)_{k}^{i} (U_\phi^m)_{j}^{k}=\delta_{j}^{i}$. Note that $(U_\phi^m)_{i j}$ and $(V_\phi^m)_{i j}$ are both symmetric. Additionally,
\begin{align}
&\big(V_\phi^m\big)_{i}^{k} (g_\rho)_{k j}=\big(U_\phi^m\big)_{i j},\qquad
(g_\rho)_{i j}^{\prime}=2 \big(P_\phi^m\big)_{i k} \big(U_\phi^m\big)_{j}^{k}.\label{g'}
\end{align}

We now relate the Levi-Civita connections ${ }^{g_\rho} \nabla$ and ${ }^{g} \nabla$, and the curvature tensors ${ }^{g_\rho} R_{i j k l}$ and ${ }^{g} R_{i j k l}$.

\begin{Lemma}Let $(g_\rho,f_\rho)$ and $(g,f)$ be defined as in Theorem~$\ref{examples-theorem}$. The Levi-Civita connections of $g_\rho$ and $g$ are related by
\begin{equation}\label{connection-conformally-flat}
{ }^{g_\rho} \nabla_{i} \eta_{j}={ }^{g} \nabla_{i} \eta_{j}-\rho \big(V_\phi^m\big)_{l}^{k} \big(P_\phi^m\big)_{i; j}^{l} \eta_{k},
\end{equation} and the curvature tensors by
\begin{equation}\label{curvature-conformally-flat}
{ }^{g_\rho} R_{i j k l}={ }^{g} R_{a b k l} \big(U_\phi^m\big)_{i}^{a} \big(U_\phi^m\big)_{j}^{b}.
\end{equation}
\end{Lemma}

\begin{proof}
First note that, since ${\rm d}P_\phi^m=0$, the right hand side of equation~\eqref{connection-conformally-flat} defines a torsion-free connection.

We now show that $g_\rho$ is parallel with respect to the connection determined by the right side of equation~\eqref{connection-conformally-flat}. On differentiating the metric, we get
\begin{equation}\label{metric-differential}
^{g}\nabla_{k} (g_\rho)_{i j}=2 \rho \big(P_\phi^m\big)_{i j; k}+2 \rho^{2} \big(P_\phi^m\big)_{(i}^{l} \big(P_\phi^m\big)_{j) l; k}.
\end{equation}
Using the definitions of $U_\phi^m$ and $V_\phi^m$ and the symmetry of $(P_\phi^m)_{ij;k}$, we get
\begin{align*}
(V_\phi^m)_{l}^{m} \big(P_\phi^m\big)_{k;(i}^{l} (g_\rho)_{j) m}=\big(P_\phi^m\big)_{k;(i}^{l} \big(U_\phi^m\big)_{j) l}=\big(P_\phi^m\big)_{i j; k}+\rho \big(P_\phi^m\big)_{k;(i}^{l} \big(P_\phi^m\big)_{j) l}.
\end{align*}
Therefore
\[
{}^g\nabla_{k} (g_\rho)_{i j}-2\rho \big(V_\phi^m\big)_{l}^{m} \big(P_\phi^m\big)_{k;(i}^{l} (g_\rho)_{j) m}=0.
\]

Combining the previous two paragraphs verifies equation~\eqref{connection-conformally-flat}.

Now set
\[
\big(D_\phi^m\big)_{j k}^{i}:=\rho \big(V_\phi^m\big)_{l}^{i} \big(P_\phi^m\big)_{j; k}^{l},
\]
so that $(D_\phi^m)_{j k}^{i}$ is the difference of ${}^{g}\nabla$ and ${}^{g_\rho}\nabla$. The difference of the curvature tensors of the connections is given in terms of $D_\phi^m$ by
\[
{ }^{g_\rho} R_{m j k l} (g_\rho)^{i m}-{ }^{g} R_{j k l}^{i}=2 \big(D_\phi^m\big)_{j[l; k]}^{i}+2 \big(D_\phi^m\big)_{j[l}^{c} \big(D_\phi^m\big)_{k] c}^{i};
\]
see \cite[p.~70]{FeffermanGraham2012}.
We compute that
\begin{align*}
\big(D_\phi^m\big)_{j l; k}^{i} &=\rho\bigl(\big(V_\phi^m\big)_{a; k}^{i} \big(P_\phi^m\big)_{j; l}^{a}+\big(V_\phi^m\big)_{a}^{i} \big(P_\phi^m\big)_{j; l k}^{a}\bigr)
\\[.5mm]
&=-\rho \big(V_\phi^m\big)_{b}^{i} \big(U_\phi^m\big)_{c; k}^{b} \big(V_\phi^m\big)_{a}^{c} \big(P_\phi^m\big)_{j; l}^{a}+\rho \big(V_\phi^m\big)_{a}^{i} \big(P_\phi^m\big)_{j; l k}^{a}
\\[.5mm]
&=-\rho^{2} \big(V_\phi^m\big)_{b}^{i} \big(P_\phi^m\big)_{c; k}^{b} \big(V_\phi^m\big)_{a}^{c} \big(P_\phi^m\big)_{j; l}^{a}+\rho \big(V_\phi^m\big)_{a}^{i} \big(P_\phi^m\big)_{j; l k}^{a}
\\[.5mm]
&=-\big(D_\phi^m\big)_{c k}^{i} \big(D_\phi^m\big)_{j l}^{c}+\rho \big(V_\phi^m\big)_{a}^{i} \big(P_\phi^m\big)_{j; l k}^{a}.
\end{align*}
Therefore,
\begin{align*}
{ }^{g_\rho} R_{m j k l} (g_\rho)^{i m} &={ }^{g} R_{j k l}^{i}+2 \rho \big(V_\phi^m\big)^{i a} \big(P_\phi^m\big)_{a j;[l k]}
\\[.5mm]
&={ }^{g} R_{j k l}^{i}+\rho \big(V_\phi^m\big)^{i a}\bigl({ }^{g} R_{a l k}^{b} \big(P_\phi^m\big)_{b j}+{ }^{g} R_{j l k}^{b} \big(P_\phi^m\big)_{a b}\bigr).
\end{align*}
Since $(V_\phi^m)_{i}^{k} (g_\rho)_{k j}=(U_\phi^m)_{i j}$, we conclude that
\begin{align*}
{ }^{g_\rho} R_{i j k l} &={ }^{g} R_{j k l}^{b} (g_\rho)_{b i}+\rho \big(U_\phi^m\big)_{i}^{a}\bigl({ }^{g} R_{a l k}^{b} \big(P_\phi^m\big)_{b j}+{ }^{g} R_{j l k}^{b} \big(P_\phi^m\big)_{a b}\bigr)
\\[.5mm]
&={ }^{g} R_{j k l}^{b} \big(U_\phi^m\big)_{b a} \big(U_\phi^m\big)_{i}^{a}+\rho \big(U_\phi^m\big)^{a}_{i}\bigl({ }^{g} R_{a l k}^{b} \big(P_\phi^m\big)_{b j}+{ }^{g} R_{j l k}^{b} \big(P_\phi^m\big)_{a b}\bigr)
\\[.5mm]
&=\bigl[{ }^{g} R_{a b k l}\bigl(\delta_{j}^{b}+\rho \big(P_\phi^m\big)_{j}^{b}\bigr)\bigr] \big(U_\phi^m\big)_{i}^{a} \\[.5mm]
&={ }^{g} R_{a b k l} \big(U_\phi^m\big)_{j}^{b} \big(U_\phi^m\big)^{a}{ }_{i}.\tag*{\qed}
\end{align*}
\renewcommand{\qed}{}
\end{proof}

The metric $\widetilde{g}$ in Theorem~$\ref{examples-theorem}$\ref{Theorem 1.2-2} is in fact flat.
\begin{Proposition}
\label{conformally-flat-theorem}
Let $(g_\rho,f_\rho)$ and $(g,f)$ be as in Theorem~$\ref{examples-theorem}$. Then $\widetilde{g}$ is flat.
\end{Proposition}

\begin{proof}
The curvature tensor of $\widetilde{g}=2\rho\,\mathrm{d}\rho\, \mathrm{d}t+t^2 g_\rho+2t\,{\rm d}\rho\, \mathrm{d}t$ is \cite[equation (6.1)]{FeffermanGraham2012}
\begin{gather*}
\widetilde{R}_{IJK0} = 0,
\\
\widetilde{R}_{i j k l} =t^{2}\bigl[{}^{g_\rho}R_{i j k l}+(g_\rho)_{i[l} (g_\rho)_{k] j}^{\prime}+(g_\rho)_{j[k} (g_\rho)_{l] i}^{\prime}-\rho (g_\rho)_{i[l}^{\prime} (g_\rho)_{k] j}^{\prime}\bigr],
\\
\widetilde{R}_{\infty j k l} =\frac{1}{2} t^{2}\bigl[\nabla_{l} (g_\rho)_{j k}^{\prime}-\nabla_{k} (g_\rho)_{j l}^{\prime}\bigr],
\\
\widetilde{R}_{\infty j k \infty} =\frac{1}{2} t^{2}\biggl[(g_\rho)_{j k}^{\prime \prime}-\frac{1}{2} (g_\rho)^{p q} (g_\rho)_{j p}^{\prime} (g_\rho)_{k q}^{\prime}\biggr].
\end{gather*}

Hence, $\widetilde{g}$ is flat if and only if
\begin{subequations}
\begin{gather}
0={ }^{g_\rho} R_{i j k l}+(g_\rho)_{i[l} (g_\rho)_{k] j}^{\prime}+(g_\rho)_{j[k} (g_\rho)_{l] i}^{\prime}-\rho (g_\rho)_{i[l}^{\prime} (g_\rho)_{k] j}^{\prime},
\label{curvature1}
\\
\label{curvature2}
0={ }^{g_\rho} \nabla_{k} (g_\rho)_{i j}^{\prime}-{ }^{g_\rho} \nabla_{j} (g_\rho)_{i k}^{\prime},
\\
\label{curvature3}
0=(g_\rho)_{i j}^{\prime \prime}-\frac{1}{2} (g_\rho)^{p q} (g_\rho)_{i p}^{\prime} (g_\rho)_{j q}^{\prime}.
\end{gather}
\end{subequations}

Equation~\eqref{curvature3} follows directly from the definition of $g_\rho$.

We next verify equation~\eqref{curvature2}.
Differentiating equation~\eqref{metric-differential} with respect to $\rho$ yields
\begin{equation}
\label{metric-double-differential}
^{g}\nabla_{k} (g_\rho)_{i j}^{\prime}=2 \big(P_\phi^m\big)_{i j; k}+4 \rho \big(P_\phi^m\big)_{(i}^{l} \big(P_\phi^m\big)_{j) l; k}.
\end{equation}
Using equation~\eqref{g'} gives
\begin{align}
\label{P^2}
\big(V_\phi^m\big)_{l}^{s} \big(P_\phi^m\big)_{k;(i}^{l} (g_\rho)_{j) s}^{\prime}=2 \big(P_\phi^m\big)_{k;(i}^{l} \big(P_\phi^m\big)_{j) l}.
\end{align}
Combining equations~\eqref{connection-conformally-flat}, \eqref{metric-double-differential} and~\eqref{P^2} yields{\samepage
\[
{ }^{g_\rho} \nabla_{k} (g_\rho)_{i j}^{\prime}=2 \big(P_\phi^m\big)_{i j; k}.
\]
Since ${\rm d}P_\phi^m=0$, we conclude that equation~\eqref{curvature2} holds.}

We now verify equation~\eqref{curvature1}. Note that
\[
g_\rho-\frac{1}{2} \rho g_\rho^{\prime}=U_\phi^m.
\]
Thus
\begin{gather*}
(g_\rho)_{i[l} (g_\rho)_{k] j}^{\prime}+(g_\rho)_{j[k} (g_\rho)_{l] i}^{\prime}-\rho (g_\rho)_{i[l}^{\prime} (g_\rho)_{k] j}^{\prime}=\big(U_\phi^m\big)_{i[l} (g_\rho)_{k] j}^{\prime}+\big(U_\phi^m\big)_{j[k} (g_\rho)_{l] i}^{\prime}.
\end{gather*}
Using equation~\eqref{g'}, we deduce that
\begin{gather*}
(g_\rho)_{i[l} (g_\rho)_{k] j}^{\prime}+(g_\rho)_{j[k} (g_\rho)_{l] i}^{\prime}-\rho (g_\rho)_{i[l}^{\prime} (g_\rho)_{k] j}^{\prime}
=2 \big(U_\phi^m\big)_{i}^{a} \big(U_\phi^m\big)_{j}{ }^{b}\bigl(g_{a[l} \big(P_\phi^m\big)_{k] b}+g_{b[k} \big(P_\phi^m\big)_{l] a}\bigr).
\end{gather*}
Since $A_\phi^m=0$, it holds that
\begin{equation}
\label{base-curvature}
{ }^{g} R_{a b k l}=-2\big(g_{a[l} \big(P_\phi^m\big)_{k] b}+g_{b[k} \big(P_\phi^m\big)_{l] a}\big).
\end{equation}
Combining equations~\eqref{curvature-conformally-flat}, \eqref{base-curvature} and~\eqref{christoffel-weighted-normal-straight} yields equation~\eqref{curvature1}.
\end{proof}

\begin{proof}[Proof of Theorem~\ref{examples-theorem}\ref{Theorem 1.2-2}]
 Proposition~\ref{conformally-flat-theorem} implies that $\widetilde{\mathrm{Ric}_\phi^m}=-\frac{m}{\widetilde{f}}\widetilde{\nabla}^2 \widetilde{f}$.

We now prove that $\widetilde{\nabla}^2 \widetilde{f}=0$. Recall \cite[equation~(3.16)]{FeffermanGraham2012} that the Christoffel symbols of $\widetilde{g}=2\rho\, \mathrm{d}t^2+t^2 g_\rho+2\rho\, \mathrm{d}t^2$ are
\begin{gather}
\widetilde{\Gamma}_{I J}^{0}=\begin{pmatrix}
0 & 0 & 0 \\
0 & -\frac{1}{2} t (g_\rho)_{i j}^{\prime} & 0 \\
0 & 0 & 0
\end{pmatrix}\!, \qquad
\widetilde{\Gamma}_{I J}^{k}=\begin{pmatrix}
0 & t^{-1} \delta_{j}^{k} & 0 \\
t^{-1} \delta_{i}^{k} & \Gamma_{i j}^{k} & \frac{1}{2} (g_\rho)^{k l} (g_\rho)_{i l}^{\prime} \\
0 & \frac{1}{2} (g_\rho)^{k l} (g_\rho)_{j l}^{\prime} & 0
\end{pmatrix}\!, \nonumber
\\
\widetilde{\Gamma}_{I J}^{\infty}=\begin{pmatrix}
0 & 0 & t^{-1} \\
0 & -(g_\rho)_{i j}+\rho (g_\rho)_{i j}^{\prime} & 0 \\
t^{-1} & 0 & 0
\end{pmatrix}\!.
\label{christoffel-weighted-normal-straight}
\end{gather}
Direct computation immediately yields
\[
\big(\widetilde{\nabla}^2 \widetilde{f}\big)_{0 I}=\big(\widetilde{\nabla}^2 \widetilde{f}\big)_{\infty\infty}=0.
\]

We next prove that $(\widetilde{\nabla}^2\widetilde{f})_{i\infty}=0$. Observe that, by Lemma~\ref{locally-conformally-flat-implications},
\begin{equation}\label{fyphim}
\frac{1}{m}\partial_i \big(fY_\phi^m\big)=\big(P_\phi^m\big)_i^k \partial_k f,
\end{equation}
and, by equation~\eqref{g'},
\[
\frac{1}{2}g^{kl}g'_{li}=\big(V_\phi^m\big)^k_s \big(P_\phi^m\big)^s_i.
\]
We then compute that{\samepage
\begin{gather*}
\begin{split}
\frac{1}{t}\big(\widetilde{\nabla}^2 \widetilde{f}\big)_{i\infty}&=\partial_i \biggl(\frac{Y_\phi^m}{m}f\biggr)-\widetilde{\Gamma}_{i\infty}^k \partial_k\biggl(f+\rho \frac{Y_\phi^m}{m}f\biggr)\\
&=\big(P_\phi^m\big)^k_i \partial_k f-\big(V_\phi^m\big)^k_s\big(P_\phi^m\big)^s_i\partial_k f-\rho\big(V_\phi^m\big)^k_s \big(P_\phi^m\big)^s_i \big(P_\phi^m\big)^l_k \partial_l f.
\end{split}
\end{gather*}
The conclusion follows from the identity $-\rho (P_\phi^m)^l_k=\delta^l_k-(U_\phi^m)^l_k$.}

We now prove that $(\widetilde{\nabla}^2\widetilde{f})_{ij}=0$. Note that
\[
\frac{1}{t}\widetilde{\nabla}_{ij}^2 f={}^{g_\rho}\nabla_{ij}^2\biggl(f+\frac{f}{m}\rho Y_\phi^m\biggr)-\frac{1}{t}\Gamma_{ij}^0\biggl(f+\frac{f}{m}\rho Y_\phi^m\biggr)-\frac{1}{m}\Gamma_{ij}^\infty fY_\phi^m.
\]
Now applying equations~\eqref{g'} and~\eqref{connection-conformally-flat} yields
\begin{gather}
\frac{1}{t}\big(\widetilde{\nabla}^2\widetilde{f}\big)_{ij}= {}^g\nabla^2_{ij}\biggl(f+\frac{f}{m}Y_\phi^m\rho\biggr)
-\rho \big(V_\phi^m\big)_l^k\big(P_\phi^m\big)^l_{i;j}\partial_k\biggl(f+\frac{f}{m}Y_\phi^m\rho\biggr)\nonumber
\\ \hphantom{\frac{1}{t}\big(\widetilde{\nabla}^2\widetilde{f}\big)_{ij}=}
{}+f\biggl(\big(P_\phi^m\big)_{ik}+\frac{1}{m}Y_\phi^m g_{ik}\biggr)\big(U_\phi^m\big)_j^k.\label{7.14}
\end{gather}
Equation~\eqref{fyphim} implies that
\begin{align}
\label{7.15}
{}^g\nabla^2_{ij}\big(fY_\phi^m\big)&=m\big[\big(P_\phi^m\big)^k_{i;j}\partial_k f+\big(P_\phi^m\big)_i^k \nabla^2_{jk}f\big].
\end{align}
Also, Lemma~\ref{locally-conformally-flat-implications} implies that
{\samepage\begin{align}
\label{7.16}
{}^g\nabla^2_{ij}f&=-f\big(P_\phi^m\big)_{ij}-\frac{f}{m}Y_\phi^mg_{ij}.
\end{align}
Equations~\eqref{7.14}, \eqref{7.15} and~\eqref{7.16}, put together, yield $(\widetilde{\nabla}^2\widetilde{f})_{ij}=0$. Hence $\widetilde{R}_{IJ}=0$.}

We now prove that $\widetilde{F_\phi^m}=0$. Since $(\widetilde{\nabla}^2\widetilde{f})_{IJ}=0$, it holds that $\widetilde{\Delta}\widetilde{f}=0$. For $m=1$, we have $\widetilde{F_\phi^m}=0$. Suppose now that $m\neq 1$. Since $(\widetilde{\nabla}^2 \widetilde{f})_{I\infty}=0$, it holds that $\partial_\rho(|\widetilde{\nabla}\widetilde{f}|^2-\mu)=0$. Hence, it suffices to show that ${(|\widetilde{\nabla}\widetilde{f}|^2-\mu)}\big|_{\rho=0}=0$. A direct computation yields
\[
\big(\big|\widetilde{\nabla} \widetilde{f}\big|^{2}-\mu\big)\big|_{\rho=0}=2\frac{f^2}{m}Y_\phi^m+|\nabla f|^2-\mu.
\]
On the one hand, it holds that~\cite[Lemma 3.1]{Case2014s}
\begin{equation}
\label{Yphim}
\frac{f^2}{m} Y_{\phi}^{m}=-\frac{1}{d+m-2}\bigl[(m-1)\big(|\nabla f|^2-\mu\big)+f\Delta f+f^2 J_\phi^m\bigr].
\end{equation}
On the other hand, the trace of equation~\eqref{second-condition-locally-conformally-flat} gives
\begin{equation}
\label{trace-of-2.2b}
f\Delta f+f^2 J_\phi^m+ \frac{d-m}{m}f^2 Y_\phi^m=0.
\end{equation}
Combining equation~\eqref{Yphim} with equation~\eqref{trace-of-2.2b} yields $\widetilde{F_\phi^m}=0$.

As both $\widetilde{\mathrm{Ric}_\phi^m}_{IJ}=0$ and $\widetilde{F_\phi^m}=0$, we conclude that $(\widetilde{\mathcal{G}},\widetilde{g},\widetilde{f},m,\mu)$ is a weighted ambient metric for $(M^d,g,f,m,\mu)$.
\end{proof}

\subsection{The Gover--Leitner conditions}

Let $(M^d,g,1,m,\mu)$, $d\geq 3$ and $m<\infty$, be a smooth metric measure space such that
\begin{equation}
\label{gover-leitner-ric}
\mathrm{Ric}_\phi^m=-(d-1)\mu g.
\end{equation} Then $F_\phi^m=-(m-1)\mu$.

Set $\lambda=-\mu/2$ and $g_\rho=(1+\lambda\rho)^2 g$ and $f_\rho=1-\lambda\rho$. Define $(\widetilde{g},\widetilde{f})$ as in equation~\eqref{weighted-straight-normal-form}. It~follows immediately from equation~\eqref{secondorderricequation} that $\widetilde{\mathrm{Ric}}_\phi^m=0$ and $\widetilde{F}_\phi^m=0$. Hence, $(\widetilde{G},\widetilde{g},\widetilde{f},m,\mu)$ is the weighted ambient metric of $(M,g,1,m,\mu)$. Moreover, equations~\eqref{fphim-equation} and~\eqref{gover-leitner-ric} imply that
\begin{align*}
&P_\phi^m=\lambda g,\qquad
Y_\phi^m=-m\lambda.
\end{align*}
This verifies Theorem~\ref{examples-theorem}\ref{Theorem 1.2-3}.

\subsection*{Acknowledgements}

We thank the anonymous referees for their valuable comments. JSC was supported by the Simons Foundation (Grant \#524601).

\pdfbookmark[1]{References}{ref}
\LastPageEnding

\end{document}